\documentclass[a4paper,twoside,11pt,reqno]{article}

\usepackage[plain]{fullpage} %full pages!
\usepackage{amsmath}
\usepackage{amssymb}
\usepackage{latexsym}
\usepackage{amsxtra}
\usepackage{amscd}
\usepackage{amsthm}
\usepackage{mathabx}
\input xy
\xyoption{all}
\entrymodifiers={+!!<0pt,\fontdimen22\textfont2>}

\usepackage{graphics}
\usepackage{epic}

\usepackage{mathrsfs, amsfonts, latexsym}%, dsfont}

\usepackage{units}

\usepackage{color}

{%\theorembodyfont{\slshape}

\newtheoremstyle{thmstyle}
 {6pt}% measure of space to leave above the theorem. E.g.: 3pt
 {6pt}% measure of space to leave below the theorem. E.g.: 3pt
 {\slshape}% name of font to use in the body of the theorem
 {}% measure of space to indent
 {\bf}% name of head font
 {.}% punctuation between head and body
 {5pt}% space after theorem head
 {}% Manually specify head

\theoremstyle{thmstyle} 
\swapnumbers 
\newtheorem{thm}{Theorem}[section]
\newtheorem{cor}[thm]{Corollary}
\newtheorem{lem}[thm]{Lemma}
\newtheorem{prop}[thm]{Proposition}

\newtheorem{defn}[thm]{Definition}

}

{%\theorembodyfont{\rmfamily}

\newtheoremstyle{defstyle}
 {6pt}% measure of space to leave above the theorem. E.g.: 3pt
 {6pt}% measure of space to leave below the theorem. E.g.: 3pt
 {}% name of font to use in the body of the theorem
 {}% measure of space to indent
 {\bf}% name of head font
 {.}% punctuation between head and body
 {5pt}% space after theorem head
 {}% Manually specify head
 
\theoremstyle{defstyle}
\swapnumbers 

\newtheorem{rem}[thm]{Remark}

}
\let\oldproofname=\proofname
\renewcommand{\proofname}{\rm\bf{\oldproofname}}

{\newtheoremstyle{sectstyle}
  {6pt}% measure of space to leave above the theorem. E.g.: 3pt
  {6pt}% measure of space to leave below the theorem. E.g.: 3pt
  {}% name of font to use in the body of the theorem
  {}% measure of space to indent
  {\bf }% name of head font        % Hier Klammer auf
  {}% punctuation between head and body % Hier Klammer zu
  {5pt}% space after theorem head
  {}% Manually specify head

\theoremstyle{sectstyle}

\newtheorem{sect}[thm]{} 
}

\renewcommand{\em}{\sl}

\newcommand{\Endproof}{\hspace*{\fill} $\Box$ \vspace{1ex} \noindent }

\makeatletter
\renewcommand{\subsection}{\@startsection{subsection}{2}%
{\z@}{-3.25ex plus -1ex minus-.2ex}{-1em}{\bf}} \makeatother

\newcommand{\LeftEqNo}{\let\veqno\leqno}

\numberwithin{equation}{section}%{subsection}
\numberwithin{thm}{section}%{subsection}
\theoremstyle{plain}

\title{Sections, Homotopy Rational Points and Reductions of Curves}
\author{Johannes Schmidt\thanks{Supported by DFG-Forschergruppe 1920 "Symmetrie, Geometrie und Arithmetik", Heidelberg--Darmstadt}}
\date{\vspace{-5ex}}

\begin{document}

\maketitle

\begin{quotation} 
\noindent \small {\bf Abstract.}
We study unramified sections of the fundamental group sequence of geometrically connected smooth projective curves of genus $\geq 2$ over $p$-adic fields together with an integral model.
We are particularly interested in the induced ``specialized'' sections of the special fibre and how they relate to homotopy rational points over the residue field.
Under mild assumptions, such a specialized section induces a unique homotopy rational point of the special fibre that is compatible with the original section of the generic fibre in cohomological settings.
We give two applications of such ``specialized'' homotopy rational points around the $\ell$-adic cycle class of a section.
\end{quotation}
\begin{quotation} 
%\\
%Keywords and Phrases: Homotopy Rational Points, Section Conjecture.
\end{quotation}

\section*{Introduction}

\subsection*{Sections and homotopy rational points.}
The profinite homotopy type of the spectrum of a field $K$ is weakly equivalent to the classification space $BG_K$ of the absolute Galois group of $K$.
Thus, any $K$-variety $X$ comes equipped with a morphism $\widehat{\rm Et}(X) \rightarrow \widehat{\rm Et}(K) \simeq BG_K$ of profinite homotopy types and each $K$-rational point of $X$ induces a splitting of this morphism.
We call any such splitting $BG_K \rightarrow \widehat{\rm Et}(X)$ in the homotopy category of simplicial profinite sets over $BG_K$ (cf.\ \cite{Quick08}) a \textbf{homotopy rational point}.
See Sect.~\ref{sect: ho rat points} for a short summary of some basic facts on homotopy rational points.

Say, $X$ is a $K(\pi,1)$-space, i.e., the homotopy type of $X$ is weakly equivalent to the classification space $B\pi_1(X)$ of its fundamental group.
This is the case e.g.\ for $X/K$ any smooth curve except for rational projective curves (see \cite{ASchmidt96} Prop.\ 15). 
In this $K(\pi,1)$-case, a homotopy rational point is equivalent (up to inner automorphisms) to a \textbf{section} of (the quotient map in) the exact sequence
\begin{equation*}\LeftEqNo\tag*{$\pi_1(X/K)$:}
 \xymatrix{
  \boldsymbol{1} \ar[r] &
  \pi_1(X\otimes_KK^{\rm s},\bar{x}) \ar[r] &
  \pi_1(X,\bar{x}) \ar[r] &
  G_K \ar[r] &
  \boldsymbol{1}
 }.
\end{equation*}

\subsection*{Specialized sections and homotopy rational points.}
Let us fix a $p$-adic field $k$ with valuation ring $\mathfrak{o}$ and residue field $\mathbb{F}$.
Further, let $X/k$ be a geometrically connected smooth projective curve of genus $\geq 2$ together with a proper flat model $\mathfrak{X}/\mathfrak{o}$ with reduced special fibre $Y = (\mathfrak{X}\otimes_{\mathfrak{o}}\mathbb{F})_{\rm red}/\mathbb{F}$.
Let $s$ be a section of the fundamental group sequence $\pi_1(X/k)$.
We say that $s$ \textbf{specializes} to a section $\bar{s}$ of $\pi_1(Y/\mathbb{F})$, if $s$ and $\bar{s}$ are compatible via the specialization map of fundamental groups:
\begin{equation*}
 \xymatrix{
  G_k \ar[r]^-s \ar[d]^-{\rm can.} &
  \pi_1(X) \ar[d]^-{\rm sp}
 \\
  G_{\mathbb{F}} \ar[r]^-{\bar{s}} &
  \pi_1(Y)
 }
\end{equation*}
(see \cite{Stix12} Ch.\ 8).
In fact, any section $s$ specializes at least to a section $\bar{s}_\ell$ of the geometrically pro-$\ell$ completed sequence $\pi_1(Y/\mathbb{F})$ (this follows from \cite{Stix12} Prop.\ 91 - see Lem.\ \ref{lem: unramified on geometric l completion}, below).
Similarly, treating $s$ as a homotopy rational point $BG_k \rightarrow X$, we say that $s$ \textbf{specializes} to a homotopy rational point $\bar{s}$ of $Y/\mathbb{F}$, if $s$ and $\bar{s}$ are compatible via the specialization morphism of homotopy types:
\begin{equation*}
 \xymatrix{
  BG_k \ar[r]^-s \ar[d]^-{\rm can.} &
  \widehat{\rm Et}(X) \ar[d]^-{\rm sp}
 \\
  BG_{\mathbb{F}} \ar[r]^-{\bar{s}} &
  \widehat{\rm Et}(Y)
 }
\end{equation*}
(see Def.\ \ref{def: specialized ho rat pt}, below).
In case of bad reduction, it is no longer clear if the special fibre $Y$ is still a $K(\pi,1)$, so the immediate question arising is,
\begin{center}
 (\dag)~ \parbox[t]{14cm}{\textit{If a section $s$ of $\pi_1(X/k)$ specializes to a section of a sufficiently nice (but not necessarily $K(\pi,1)$-) reduction $Y/\mathbb{F}$, does such a specialized section induce a specialized homotopy rational point of $Y/\mathbb{F}$ ?}}
\end{center}
Say, the section $s$ of the fundamental group sequence $\pi_1(X/k)$ specializes to a section $\bar{s}$ of $\pi_1(Y/\mathbb{F})$.
Suppose all the points in the normalization $\pi: \tilde{Y} \rightarrow Y$ lying above singular points of $Y$ contained in rational components are $\mathbb{F}$-rational.
Studying the homotopy type of $Y/\mathbb{F}$ (see Thm.\ \ref{thm: homotopy type of a curve}, below), at least, a possible candidate for a specialized homotopy rational point of $s$ giving back the specialized section $\bar{s}$ does exist (see Cor.\ \ref{cor: sections vs ho rational points}).
However, the compatibility of $s$ and this candidate for a specialized homotopy rational point (in abuse of notation also denoted by) $\bar{s}$ via the specialization morphism of homotopy types is not clear.
We will prove this compatibility at least in the sufficiently additive setting of cohomology cochains:
Let $\Lambda$ be a discrete torsion $G_{\mathbb{F}}$-module.
The homotopy rational point $s$ resp.\ the candidate $\bar{s}$ induces a spitting $s^*$ resp.\ $\bar{s}^*$ of the canonical morphism $\mathbb{R}\Gamma(k_{\rm \acute{e}t},\Lambda) \rightarrow \mathbb{R}\Gamma(X_{\rm \acute{e}t},\Lambda)$ resp.\ $\mathbb{R}\Gamma(\mathbb{F}_{\rm \acute{e}t},\Lambda) \rightarrow \mathbb{R}\Gamma(Y_{\rm \acute{e}t},\Lambda)$ in the derived category.
It turns out that these splittings are compatible via the specialization morphism of cohomology cochains (for the exact statement, see Thm.\ \ref{thm: specialized homotopy rational point}, below): 

\bigskip\noindent {\bf Theorem A.} {\it
Let $X/k$ be a geometrically connected smooth projective curve of genus $\geq 2$ over a $p$-adic field $k$ together with a regular, proper, flat model $\mathfrak{X}/\mathfrak{o}$ with reduced special fibre $Y=(\mathfrak{X}\otimes_{\mathfrak{o}}\mathbb{F})_{\rm red}/\mathbb{F}$.
Suppose all the points in the normalization $\pi: \tilde{Y} \rightarrow Y$ lying above singular points on rational components of $Y$ are $\mathbb{F}$-rational.
Let $s$ be a section of $\pi_1(X/k)$ specializing to a section of $\pi_1(Y/\mathbb{F})$ and let $\Lambda$ be a constructible $G_{\mathbb{F}}$-module.
Then there is a homotopy rational point $\bar{s}$ of $Y/\mathbb{F}$ inducing the specialized section on fundamental groups and a commutative diagram of cohomology cochains in the derived category $\mathcal{D}^+(\underline{\rm Ab})$:
\begin{equation*}
 \xymatrix{
  \mathbb{R}\Gamma(X_{\rm \acute{e}t},\Lambda) \ar[r]^-{s^*} &
  \mathbb{R}\Gamma(k_{\rm \acute{e}t},\Lambda) \phantom{.}
 \\
  \mathbb{R}\Gamma(Y_{\rm \acute{e}t},\Lambda) \ar[r]^-{\bar{s}^*} \ar[u]_-{{\rm sp}^*} &
  \mathbb{R}\Gamma(\mathbb{F}_{\rm \acute{e}t},\Lambda). \ar[u]
 }
\end{equation*}
{\it}}

\subsection*{Applications.}
We will give two applications around the cycle class ${\rm cl}_s$ of a section (see \cite{Stix12} Sect. 6.1, \cite{EsnaultWittenberg09} or Rem.~\ref{rem: cycle class}, below).
The first is an application of specialized sections in general\footnote{In fact, it was the original motivation for question (\dag).} and the second an application of Thm.~A itself.

Let us shortly describe the first application:
Based around a pro-$\ell$ specialization result, we give an independent proof to the following algebraicity result of Esnault and Wittenberg:

\bigskip\noindent {\bf Proposition B.} (\cite{EsnaultWittenberg09} Cor.\ 3.4) {\it
Let $X/k$ be a geometrically connected smooth projective curve of genus $\geq 2$ over a $p$-adic field $k$, admitting a section $s$ of $\pi_1(X/k)$.
Then the $\ell$-adic cycle class ${\rm cl}_s$ of $s$ lies in the image of the Chern class map $\hat{c}_1:{\rm Pic}(X)\otimes \mathbb{Z}_\ell \rightarrow {\rm H}^2(X,\mathbb{Z}_\ell(1))$ for each prime $\ell \neq p$.
{\it}}

\bigskip
Using Tate-Lichtenbaum duality, this prime-to-$p$ algebraicity is equivalent to $s^*$ trivializing the first Chern class map $\hat{c}_1$ in ${\rm H}^2(G_k,\mathbb{Z}_\ell(1))$ for $\ell \neq p$.
This is precisely the statement we will prove in Prop.~\ref{prop: The l-adic chern class mapping is killed by a section}, below.
In fact, the pullback along $s$ of certain Chern classes will turn out to be exactly the obstructions for the compatibility claimed in Thm.~A.
Thus, we have to avoid the usage of Thm.~A in the proof of Prop.~B resp.\ the equivalent Prop.~\ref{prop: The l-adic chern class mapping is killed by a section} and modify $X/k$ until it admits a $K(\pi,1)$-model over $\mathfrak{o}$.

Our second application is a direct consequence of (a geometrically pro-$\ell$ completed variant of) Thm.\ A:
In \cite{EsnaultWittenberg09} Rem.\ 3.7 (iii) Esnault and Wittenberg raised the question whether the $\ell$-adic cycle class ${\rm cl}_s$ of a section $s$ of $\pi_1(X/k)$ admits a canonical lift to ${\rm H}^2(\mathfrak{X},\mathbb{Z}_\ell(1))$.
Using the (unconditional) geometrically pro-$\ell$ completed specialization $\bar{s}_\ell$ of $s$, we construct a canonical cycle class ${\rm cl}_s^{\mathfrak{X}}$ in ${\rm H}^2(\mathfrak{X},\mathbb{Z}_\ell(1))$ for $\mathfrak{X}/\mathfrak{o}$ a regular model satisfying the assumptions of Thm.\ A.
A (geometrically pro-$\ell$ completed) variant of Thm.~A shows that ${\rm cl}_s^{\mathfrak{X}}$ is indeed a lift of ${\rm cl}_s$ (see Prop.\ \ref{prop: canonical lift of the cycle class to the model}, below):

\bigskip\noindent {\bf Proposition C.} {\it
Let $X/k$ be a geometrically connected smooth projective curve of genus $\geq 2$ over a $p$-adic field $k$ together with a regular, proper, flat model $\mathfrak{X}/\mathcal{O}$ satisfying the assumptions of Thm.\ A.
Then for any $\ell \neq p$ and any section $s$ of $\pi_1(X/k)$, the induced $\ell$-adic cycle class ${\rm cl}_s$ admits a canonical lift ${\rm cl}_s^{\mathfrak{X}}$ to ${\rm H}^2(\mathfrak{X},\mathbb{Z}_\ell(1))$.
{\it}}

\subsection*{Notation.}
In the following, $k$ always denotes a $p$-adic field with valuation ring $\mathfrak{o}$ and residue field $\mathbb{F}$.
If $K$ is any field, we denote a fixed separable closure by $K^{\rm s}$ and the corresponding absolute Galois group by $G_K$.
Denote the maximal unramified subextension of $k^{\rm s}/k$ by $k^{\rm nr}/k$.
For $X$ a $K$- resp.\ $k$-variety, denote its base-change to $K^{\rm s}$ resp.\ $k^{\rm nr}$ by $X^{\rm s}$ resp.\ $X^{\rm nr}$.
Denote by $\hat{\mathcal{S}}_{(*)}$ the category of (pointed) simplicial profinite sets together with the model structures of \cite{Quick08}.
For $X$ a scheme together with a geometric point $\bar{x}$, $\pi_1(X,\bar{x})$ denotes its profinite \'{e}tale fundamental group.
Mostly we will skip the base-point in our notation.
Denote by $\widehat{\rm Et}(X)$ its profinite \'{e}tale homotopy type in $\hat{\mathcal{S}}_{(*)}$.
%If no confusion could arise, we will just write $\hat{X}$ or $X$ for $\widehat{\rm Et}(X)$.
For a diagram $\mathcal{Y}^\prime \leftarrow \mathcal{Y} \rightarrow \mathcal{Y}^{\prime\prime}$ of simplicial (profinite) sets, we write $\mathcal{Y}^\prime \vee_\mathcal{Y} \mathcal{Y}^{\prime\prime}$ for the homotopy pushout.
For a simplicial profinite set $\mathcal{Y}$ with torsion local system $\Lambda$, write $C^\bullet(\mathcal{Y},\Lambda)$ for its cohomology cochains (see \cite{Quick08} Sect.\ 2.2).
If $\mathcal{Y}$ is the homotopy type $\widehat{\rm Et}(X)$ of a scheme $X$, then $C^\bullet(\widehat{\rm Et}(X),\Lambda)$ is quasi-isomorphic to $\mathbb{R}\Gamma(X_{\rm \acute{e}t},\Lambda)$ (see \cite{Quick08} Sect.\ 3.1).
If $\mathcal{Y}$ is the classification space $BG$ of a profinite group $G$, then $C^\bullet(BG,\Lambda)$ is quasi-isomorphic to $\mathbb{R}\Gamma(G,\Lambda)$.
We will write just $C^\bullet(X,\Lambda)$ resp.\ $C^\bullet(G,\Lambda)$ in these cases.
Similarly, we will just write $C^\bullet(X\otimes_kk^\prime,\Lambda)$ for the $G$-equivariant cochains $C^\bullet(\widehat{\rm Et}(X) \times_{BG}EG,\Lambda)$ for $K^\prime/K$ a Galois extension with group $G$ and $X/K$ a $K$-variety. 
Finally, we will always use continuous \'{e}tale cohomology in the sense of \cite{Jannsen88}.

\subsection*{Acknowledgements.}
I would like to thank H\'{e}l\`{e}ne Esnault, Armin Holschbach, Gereon Quick, Alexander Schmidt and Jakob Stix for helpful comments, discussions and/or suggestions. 

\section{Preliminaries: Homotopy rational points}\label{sect: ho rat points}

\subsection*{Homotopical algebra.}
We will work in the following homotopy categories:

\begin{sect}\label{para: homotopy categories}
Let $\hat{\mathcal{S}}$ be the category of simplicial profinite sets together with the model structures of \cite{Quick08}.
For $G$ a profinite group, let $BG$ be its profinite classification space and $\hat{\mathcal{S}}\downarrow BG$ the category of simplicial profinite sets over $BG$ together with the induced model structure.
A simplicial profinite $G$-set is a simplicial profinite set together with a degreewise continuous $G$-action.
Let $\hat{\mathcal{S}}_G$ be the resulting category together with the model structure of \cite{Quick10}.
By \cite{Quick10} Cor.\ 2.11, $\mathcal{H}(\hat{\mathcal{S}}\downarrow BG)$ is Quillen equivalent to $\mathcal{H}(\hat{\mathcal{S}}_G)$ via the base change functor $\mathcal{X} \mapsto \mathcal{X} \times_{BG}EG$, where $BG$ is the profinite classification space of $G$ and $EG \rightarrow BG$ the universal covering.
Under this equivalence, maps $BG\rightarrow \mathcal{X}$ in $\mathcal{H}(\hat{\mathcal{S}}\downarrow BG)$ correspond to homotopy fixed points of $\mathcal{X} \times_{BG}EG$, i.e., maps ${\rm pt} \simeq EG \rightarrow \mathcal{X} \times_{BG}EG$ in $\mathcal{H}(\hat{\mathcal{S}}_G)$.
\end{sect}

\begin{sect}\label{para: ho fixed points}
Let $\mathcal{X}$ be a simplicial profinite $G$-set.
The set of homotopy fixed points $[EG,\mathcal{X}]_{\hat{\mathcal{S}}_G}$ is the set of connected components of Quick's homotopy fixed point space $\mathcal{X}^{hG} = S_G(EG,\mathcal{X}^\prime)$ (where $\mathcal{X}\rightarrow \mathcal{X}^\prime$ is a fibrant replacement in $\hat{\mathcal{S}}_G$), defined and studied in \cite{Quick10}.
In general, $\mathcal{X}^{hG}$ is difficult to describe.
At least, by \cite{Quick10} Thm.\ 2.16, there is a Bousfield-Kan type descent spectral sequence (with differentials in the usual ``cohomological'' directions)
\begin{equation}\label{eq: descent spectral sequence}
 E_2^{p,q} = {\rm H}^p(G,\pi_{-q}(\mathcal{X})) \Rightarrow \pi_{-(p+q)}(\mathcal{X}^{hG}).
\end{equation}
\end{sect}

Applying Bousfield and Kan's connectivity lemma \cite{BousfieldKan} Ch.\ IX 5.1 to the spectral sequence in \ref{para: ho fixed points}, one can prove: %(see \cite{JSchmidt15} Lem.~1.1 for details):

\begin{lem}\label{lem: injection of hofixed points}
Let $G$ be a profinite group of cohomological dimension $\leq n$ and $f:\mathcal{X} \rightarrow \mathcal{Y}$ an $(n+1)$-equivalence in $\hat{\mathcal{S}}_G$ (i.e., $\pi_q(f)$ is an isomorphism for all $q \leq n$ and an epimorphism for $q = n+1$).
Then $f$ induces an injection
\begin{equation*}
 [EG,\mathcal{X}]_{\hat{\mathcal{S}}_G} =
 \xymatrix{
  \pi_0(\mathcal{X}^{hG}) \ar[r] &
  \pi_0(\mathcal{Y}^{hG})
 }
 = [EG,\mathcal{Y}]_{\hat{\mathcal{S}}_G}.
\end{equation*}
\end{lem}

\proof
We may assume that $\mathcal{Y}$ is fibrant and $f$ is a fibration in $\hat{\mathcal{S}}_G$.
Further, we may assume that $\mathcal{X}^{hG}$ is non-empty.
Say, $s:E\Gamma \rightarrow \mathcal{X}$ is a model of a homotopy fixed point and let $r$ be $f \circ s$.
The fibre $\mathcal{F}_s:= \mathcal{X} \times_{\mathcal{Y}} EG$ comes equipped with a fibration into $EG$, hence is fibrant in $\hat{\mathcal{S}}_G$, too.
Taking limits (i.e., forgetting the topology in \cite{Quick10}) resp.\ simplicial mapping spaces $S_\Gamma(EG,-)$ of $\hat{\mathcal{S}}_G$ gives us a homotopy fibre sequence
\begin{equation*}
 \xymatrix{
  \lim\mathcal{F}_s \ar[r] &
  \lim\mathcal{X} \ar[r] &
  \lim\mathcal{Y}
 }
\end{equation*}
resp.\
\begin{equation*}
 \xymatrix{
  \mathcal{F}_s^{hG} \ar[r] &
  \mathcal{X}^{hG} \ar[r] &
  \mathcal{Y}^{hG}
 }
\end{equation*}
in $\underline{\rm SSets}_\ast$ (pointed by the neutral element in $G$).
By \cite{Quick13} Lem. 2.9, the limit of $f$ is an $n$-equivalence of simplicial sets. 
So, again by loc.\ cit., the first fibre sequence implies the $n$-connectedness of $\mathcal{F}_s$.
Using the second homotopy fibre sequence, we get that the map of pointed sets $(\pi_0(\mathcal{X}^{hG}),s) \rightarrow (\pi_0(\mathcal{Y}^{hG}),r)$ has kernel $\pi_0(\mathcal{F}_s^{hG})$.
Bousfield and Kan's connectivity lemma applied to the descent spectral sequence (\ref{eq: descent spectral sequence}) for $\mathcal{F}_s$ implies that this kernel is trivial, since $\mathcal{F}_s$ is $n$-connected and $G$ has cohomological dimension $\leq n$.
Varying over all the homotopy fixed points of $\mathcal{X}$, we get the result.
\Endproof

\begin{sect}\label{para: sect and vs ho rat points of classification spaces}
Let $p:\pi \twoheadrightarrow G$ be an epimorphism of profinite groups with kernel $\bar{\pi}\unlhd \pi$.
By the adjunction between the profinite groupoid $\Pi(-)$ and $B(-)$, $B\pi / BG$ is fibrant in $\hat{\mathcal{S}}\downarrow BG$.
Thus, $[BG,B\pi]_{\hat{\mathcal{S}}\downarrow BG}$ is given as ${\rm Hom}_{\hat{\mathcal{S}}\downarrow BG}(BG,B\pi)$ modulo homotopy equivalences over $BG$ with respect to the standard cylinder object $BG \otimes \Delta^1$.
Such homotopies between maps $B(s_0)$ and $B(s_1)$ for sections $s_i$ of $p$ correspond precisely to conjugation of these sections via elements of $\bar{\pi}$.
In particular, $B(-)$ resp.\ $\pi_1(-,\ast)$ give canonical identifications between $[BG,B\pi]_{\hat{\mathcal{S}}\downarrow BG}$ and the set of $\bar{\pi}$-conjugacy classes of section of $p$.
\end{sect}

\begin{sect}\label{para: sections vs ho rat points}
Let $\mathcal{X} /BG$ be a connected simplicial profinite set in $\hat{\mathcal{S}}\downarrow BG$ and assume $\pi_1(\mathcal{X}) \rightarrow G$ is an epimorphism.
Then any map $BG \rightarrow \mathcal{X}$ in $\mathcal{H}(\hat{\mathcal{S}}\downarrow BG)$ defines a $\pi_1(\mathcal{X}\times_{BG}EG)$-conjugacy class of splittings of the fundamental group sequence
\begin{equation*}
 \xymatrix{
  \boldsymbol{1} \ar[r] &
  \pi_1(\mathcal{X}\times_{BG}EG) \ar[r] &
  \pi_1(\mathcal{X}) \ar[r] &
  G \ar[r] &
  \boldsymbol{1}
 }
\end{equation*}
of $\mathcal{X} /BG$.
%Similar to the case of simplicial sets, for a simplicial profinite set $\mathcal{Y}$ and a profinite group $\pi$, homotopy classes $[\mathcal{Y},B\pi]_{\hat{\mathcal{S}}}$ correspond to functors of profinite categories modulo natural transformations between the fundamental groupoid $\Pi(\mathcal{Y})$ and $\pi$.
Conversely, if the underlying simplicial profinite set $\mathcal{X}$ of $\mathcal{X} /BG$ is a $K(\pi,1)$ (i.e., the canonical map $\mathcal{X} \rightarrow B\Pi(\mathcal{X})$ into the classification space of the profinite fundamental groupoid is a weak equivalence), it follows from \ref{para: sect and vs ho rat points of classification spaces} that sections of the above fundamental group sequence of $\mathcal{X} /BG$ modulo conjugation correspond to maps $BG \rightarrow \mathcal{X}$ in $\mathcal{H}(\hat{\mathcal{S}}\downarrow BG)$.
\end{sect}

Combining \ref{para: sections vs ho rat points} with Lem.\ \ref{lem: injection of hofixed points}, we get:

\begin{cor}\label{cor: sections vs ho rat points}
Let $G$ be a profinite group of cohomological dimension $1$ and $\mathcal{X} /BG$ a connected simplicial profinite set in $\hat{\mathcal{S}}\downarrow BG$ s.t.\ $\pi_1(\mathcal{X}) \rightarrow G$ is an epimorphism.
Assume that the canonical map $\mathcal{X}\times_{BG}EG \rightarrow B\Pi(\mathcal{X}\times_{BG}EG)$ admits a section in $\mathcal{H}(\hat{\mathcal{S}}_G)$.
Then we get a canonical identification between the set of $\pi_1(\mathcal{X}\times_{BG}EG)$-conjugacy classes of sections of $\pi_1(\mathcal{X}) \rightarrow G$ and $[BG,\mathcal{X}]_{\hat{\mathcal{S}}\downarrow BG} \simeq [EG,\mathcal{X}\times_{BG}EG]_{\hat{\mathcal{S}}_G}$.  
\end{cor}

\proof
Indeed, the canonical map $\mathcal{X}\times_{BG}EG \rightarrow B\Pi(\mathcal{X}\times_{BG}EG)$ is a $2$-equivalence, so it induces an injection on the respective sets of homotopy fixed points by Lem.\ \ref{lem: injection of hofixed points}.
Since it admits a section, the induced map on homotopy fixed points is even bijective and the claim follows from \ref{para: sections vs ho rat points}.
\Endproof

\subsection*{Homotopy rational points.}
We are mainly interested in profinite homotopy types of varieties over a field $K$:

\begin{sect}\label{para: profinite ho type}
Let $Z$ be a $K$-variety.
The profinite \'{e}tale homotopy type of the spectrum of $K$ (pointed by the choice of a separable closure $K^{\rm s}/K$) is a $K(\pi,1)$ with fundamental group $G_K$, i.e., it is weakly equivalent to the profinite classification space $BG_K$ in the pointed category $\hat{\mathcal{S}}_{\ast}$.
We define the profinite homotopy type $\widehat{\rm Et}(Z/K) \rightarrow BG_K$ of $Z/K$ as the resulting map in $\hat{\mathcal{S}}\downarrow BG_K$ induced by the structural map of $Z/K$.
Using \cite{Quick10} Thm.~3.5 and Lem.~3.3, we see that the underlying homotopy type of $\widehat{\rm Et}(Z/K)\times_{BG_K}EG_K$ together with its abstract $G_K$-action corresponds to the homotopy type $\widehat{\rm Et}(K\otimes_KK^{\rm s})$ together with the induced abstract $G_K$-action.
Similar arguments work for any Galois extension $K^\prime/K$ and $Z\otimes_KK^\prime$, too.
\end{sect}

\begin{sect}\label{para: ho rat points}
Each $K$-rational point of $Z$ defines a map $BG_K \rightarrow \widehat{\rm Et}(Z/K)$ in $\mathcal{H}(\hat{\mathcal{S}}\downarrow BG_K)$, i.e., a homotopy fixed point of $\widehat{\rm Et}(Z/K) \times_{BG_K}EG_K$.
Thus, for any simplicial profinite set $\mathcal{X}/BG_K$, we call any map $BG_K \rightarrow \mathcal{X}$ in $ \mathcal{H}(\hat{\mathcal{S}}\downarrow BG_K)$ a {\bf homotopy rational point} of $\mathcal{X}$ over $K$.
A homotopy rational point of $Z/K$ simply is a homotopy rational point of $\widehat{\rm Et}(Z/K)$.  
\end{sect}

\begin{sect}\label{para: sections vs ho rat points of varieties}
Let $Z/K$ be a geometrically connected $K$-variety.
Then any homotopy rational point of $Z/K$ gives a conjugacy class of splittings of the fundamental group sequence $\pi_1(Z/K)$.
Conversely, assume $Z$ has the $K(\pi,1)$-property, i.e., its \'{e}tale cohomology of constructible locally constant coefficients is given by the cohomology of its finite \'{e}tale site.
It follows that $\widehat{\rm Et}(Z/K)$ is a $K(\pi,1)$-space in the above sense.
By \ref{para: sections vs ho rat points}, we get a canonical identification between the set of $\pi_1(Z\otimes_KK^{\rm s})$-conjugacy classes of the fundamental group sequence $\pi_1(X/K)$ and homotopy rational resp.\ homotopy fixed points of $Z/K$ (cf.\ \cite{Quick10} Sect.~3.2).
By \cite{Stix02} Prop.\ A.4.1, this in particular is the case for $Z$ any smooth curve over $K$ except for Brauer-Severi curves.
\end{sect}

%Let us conclude our discussion of homotopy rational points with the following lemma about homotopy classes of $K$-rational points:

\begin{lem}\label{lem: ho classes of rational points}
Let $Z$ be a non-empty irreducible $K$-variety.
Then the canonical map
\begin{equation*}
  \xymatrix{
   Z(K) \ar[r] &
   [B\Gamma_K,\widehat{\rm Et}(Z/K)]_{\hat{\mathcal{S}}\downarrow B\Gamma_K}
  }
 \end{equation*}
 is trivial in the following two situations:
\begin{enumerate}
 \item If $K$ has cohomological dimension $\leq n$ and $\widehat{\rm Et}(Z/K) \times_{BG_K}EG_K$ is $n$-connected.
 \item If $K$ has characteristic $0$ and $Z$ is $\mathbb{A}^1$-chain connected (i.e., for any two $K$-rational points $z^\prime,z^{\prime\prime}$, there are finitely many $K$-morphisms $u_i: \mathbb{A}_K^1 \rightarrow Z$, $1\leq i \leq n$ with $u_1(0) = z^\prime$, $u_n(1) = z^{\prime\prime}$ and $u_i(1) = u_{i+1}(0)$). %or $Z/K$ is smooth and $\mathbb{A}^1$-connected,
\end{enumerate}
\end{lem}

\proof
To prove (i), it is enough to show that $\pi_0((\widehat{\rm Et}(Z/K)\times_{B\Gamma_K}E\Gamma_K)^{h\Gamma_K})$ is trivial.
But this follows from Bousfield and Kan's connectivity lemma applied to the descent spectral sequence (\ref{eq: descent spectral sequence}) for $\widehat{\rm Et}(Z/K)\times_{B\Gamma_K}E\Gamma_K$:
${\rm cd}(K) \leq n$ and $\widehat{\rm Et}(Z/K) \times_{BG_K}EG_K$ is $n$-connected by assumption.
Statement (ii) holds, since $\mathbb{A}_K^1$ is contractible (to $B\Gamma_K$) in characteristic $0$.
%The $\mathbb{A}^1$-connected smooth case holds, since our map from rational to homotopy-rational points factors through the map given by \'{e}tale realization of $\mathbb{A}^1$-homotopy theory (see \cite{Isaksen04} or \cite{ASchmidt11}).
\Endproof

\subsection*{Geometric pro-\boldmath{$\ell$} completions.}
Let us shortly discuss pro-$\ell$-completion in $\mathcal{H}({\hat{\mathcal{S}}}_G)$ for a strongly complete profinite group $G$.

\begin{sect}\label{para: geometric pro l completion}
In \cite{Quick12}, Quick gave an explicit construction of a pro-finite completion in $\hat{\mathcal{S}}_G$.
An analogue construction gives a pro-$\ell$-completion in $\hat{\mathcal{S}}_G$, too (see loc.~cit.~Rem.~3.3).
Let us shortly describe this construction:
By loc.~cit.~4.3, any profinite $G$-space $\mathcal{X}$ is isomorphic to a profinite $G$-space of the form $\{ \mathcal{X}_i \}_{i\in I}$ for $\mathcal{X}_i$ a finite discrete $G$-space.
Then the pro-$\ell$-completion is given as the profinite $G$-space
\begin{equation*}
 \mathcal{X}_\ell^\wedge := \{ \bar{W} \hat{\Gamma}_\ell(\mathcal{X}_i) \}_{i\in I},
\end{equation*}
where $\bar{W}(-)$ is (levelwise) the classification space and $\hat{\Gamma}_\ell(\mathcal{X}_i)$ is degreewise the pro-$\ell$ completion of the free loop group $\Gamma(\mathcal{X}_i)$ of $\mathcal{X}_i$.
Arguing directly using the (levelwise) homotopy fibre sequence
\begin{equation*}
 \xymatrix{
  \hat{\Gamma}_\ell(\mathcal{X}_i) \ar[r] &
  W \hat{\Gamma}_\ell(\mathcal{X}_i) \ar[r] &
  \bar{W}\hat{\Gamma}_\ell(\mathcal{X}_i)
 },
\end{equation*}
or comparing $\mathcal{X}_\ell^\wedge$ with the fibrant replacement in Morel's pro-$\ell$ model structure in \cite{Morel96} (see Sect.~2.1 in loc.~cit.), we get that $\pi_1(\mathcal{X}_\ell^\wedge)$ equals the pro-$\ell$ completion $\pi_1^\ell(\mathcal{X})$, $\mathcal{X}_\ell^\wedge$ has pro-$\ell$ profinite homotopy groups and the canonical map $\mathcal{X} \rightarrow \mathcal{X}_\ell^\wedge$ induces an isomorphism in $\mathcal{D}^+(\underline{\rm Mod}_G)$ on cohomology cochains $C^\bullet(-,\Lambda)$ for any finite $\ell$-torsion $G$-module $\Lambda$.
\\
For $\mathcal{X}/BG$ in $\hat{\mathcal{S}}\downarrow BG$, denote by $\mathcal{X}_{(\ell)}^\wedge / BG$ the homotopy type in $\mathcal{H}(\hat{\mathcal{S}}\downarrow BG)$ corresponding to the pro-$\ell$ completion $(\mathcal{X} \times_{BG}EG)_\ell^\wedge$ in $\mathcal{H}(\hat{\mathcal{S}}_G)$. 
\end{sect}

If $G$ itself is not a pro-$\ell$ group, $\mathcal{X}_{(\ell)}^\wedge / BG$ corresponds to a ``geometric'' pro-$\ell$ completion in the relative homotopy category $\mathcal{H}(\hat{\mathcal{S}}\downarrow BG)$.
Let us discuss the case of $B\pi \rightarrow BG$ for suitable $\pi\twoheadrightarrow G$:

\begin{sect}\label{para: geometric pro l completion of a group}
Let $p:\pi\twoheadrightarrow G$ be an epimorphism of profinite groups with kernel $\bar{\pi} \unlhd \pi$.
Assume $\bar{\pi}$ is an $\ell$-good profinite group, i.e., the pro-$\ell$ completion map $\bar{\pi} \rightarrow \hat{\bar{\pi}}_\ell$ induces isomorphisms ${\rm H}^q(\bar{\pi},\Lambda) \simeq {\rm H}^q(\hat{\bar{\pi}}_\ell,\Lambda)$ for all finite $\ell$-torsion $\hat{\bar{\pi}}_\ell$-modules $\Lambda$ and all these cohomology groups are finite.
Let $\Delta_\ell \unlhd_c \bar{\pi}$ be the kernel of the pro-$\ell$ completion of $\bar{\pi}$.
Note that it is a characteristic subgroup by the universal property of the completion.
Then we define the geometric pro-$\ell$ completion of $\pi\rightarrow G$ as $\hat{\pi}_{(\ell)}:= \pi / \Delta_\ell \rightarrow G$.
By construction, the geometric pro-$\ell$ completion $\pi \twoheadrightarrow \hat{\pi}_{(\ell)}$ sits in the following commutative diagram with exact rows
\begin{equation*}
 \xymatrix{
  \phantom{.}\boldsymbol{1} \ar[r] &
  \bar{\pi} \ar[r] \ar@{->>}[d] &
  \pi \ar[r] \ar@{->>}[d] &
  G \ar[r] \ar@{=}[d] &
  \boldsymbol{1}\phantom{.}
 \\
  \phantom{.}\boldsymbol{1} \ar[r] &
  \hat{\bar{\pi}}_\ell \ar[r] &
  \hat{\pi}_{(\ell)} \ar[r] &
  G \ar[r] &
  \boldsymbol{1}.
 }
\end{equation*}
It follows that $B\hat{\pi}_{(\ell)}\times_{BG}EG \rightarrow (B\hat{\pi}_{(\ell)}\times_{BG}EG)_\ell^\wedge$ is a weak equivalence in $\hat{\mathcal{S}}_G$ (this is a special case of the pro-$\ell$ analogue of \cite{Quick12} Thm.\ 3.14).
Further, the canonical Map $\pi \twoheadrightarrow \hat{\pi}_{(\ell)}$ induces an isomorphism $(B\pi\times_{BG}EG)_\ell^\wedge \rightarrow B\hat{\pi}_{(\ell)}\times_{BG}EG$.
In particular, $(B\pi)_{(\ell)}^\wedge = B (\hat{\pi}_{(\ell)})$.
\end{sect}

For fundamental groups of geometrically connected $K$-varieties we define:

\begin{sect}\label{para: geometric pro l completion of fundamental groups}
%Let $K$ be a field whose absolute Galois group $G_K$ is strongly complete (e.g., $K$ a finite or $p$-adic field) and let $Z/K$ be a geometrically connected $K$-variety.
Let $Z/K$ be a geometrically connected $K$-variety.
Then we write $\pi_1^\ell(Z)$ resp.\ $\pi_1^{(\ell)}(Z)$ for the pro-$\ell$ resp.\ geometric pro-$\ell$ completion of the \'{e}tale fundamental group $\pi_1(Z)$ in the sense of \ref{para: geometric pro l completion of a group}.
Denote by $\pi_1^{(\ell)}(Z/\mathbb{F})$ the geometrically pro-$\ell$ completed fundamental group sequence sequence
\begin{equation*}
 \xymatrix{
  \phantom{.}\boldsymbol{1} \ar[r] &
  \pi_1^\ell(Z\otimes_{K}K^{\rm s}) \ar[r] &
  \pi_1^{(\ell)}(Z) \ar[r] &
  G_{K} \ar[r] &
  \boldsymbol{1}.
 }
\end{equation*}
Say $G_K$ is strongly complete (e.g., $K$ a finite or $p$-adic field).
Then $\pi_1^{(\ell)}(Z)$ is the fundamental group of the geometric pro-$\ell$ completion $\widehat{\rm Et}(Z/K)_{(\ell)}^\wedge$.
\end{sect}

\section{The \'{e}tale homotopy type of a curve}\label{sect: homotopy type of a curve}

\subsection*{The \'{e}tale homotopy type of a curve}
We want to discuss the the \'{e}tale homotopy type of a curve $Z$ over a field $K$.

\begin{sect}\label{para: semi and weak normalizations}
Recall that for $Z$ a reduced $K$-variety, its semi- resp.\ weak-normalization $\pi^{\rm sn}: Z^{\rm sn} \rightarrow Z$ resp.\ $\pi^{\rm wn}: Z^{\rm wn} \rightarrow Z$
is universal among factorizations $f:Z^\prime \rightarrow Z$ of the normalization $\pi: \tilde{Z} \rightarrow Z$ with $f$ birational, bijective on points and inducing trivial resp.\ purely inseparable extensions on the residue fields of all points (cf.\ \cite{Kollar96} Sect.\ I.7.2).
In particular, for $Z$ a curve over a perfect field, the semi- and weak-normalization coincide.
Further, the curve $Z$ itself is semi-normal (i.e., agrees with its semi-normalization), if and only if geometrically it has at most ordinary multiple points as singularities (use \cite{Kollar96} Sect.\ I 7.2.2.1). %see also loc.~cit.~Prop.~7.2.6.
\end{sect}

\begin{sect}\label{para: dual graph}
Let $Z$ be a reduced curve over $K$ and assume either $Z$ to be semi-normal or $K$ to be perfect. 
Denote by $\Gamma(Z)$ the following bipartite graph:
The two sets of nodes are the irreducible components resp.\ the singular points of $Z^{\rm sn}$ and a node corresponding to a component $Z_i$ is joined to a node corresponding to a singularity $z$ by multiplicity-of-$z$-in-$Z_i$-many edges. 
Let $\Gamma.(Z)$ in $\hat{\mathcal{S}}$ be the profinite completion of the canonical realization of $\Gamma(Z)$ as a simplicial set.
We will call both $\Gamma(Z)$ and $\Gamma.(Z)$ the \textbf{dual-(bipartite)graph} of $Z$.
Obviously, $\Gamma.(Z)$ has the homotopy type of the profinite completion of a bouquet of $S^1$'s, i.e., is isomorphic to $B F_r$ for $F_r$ the free profinite group of $r$ generators and $r$ the number of loops in $\Gamma(Z)$.
If $Z$ is even semi-stable, then $\Gamma(Z)$ is just the barycentric subdivision of the regular dual-graph of the semi-stable curve $Z$.
It follows from the universal property of weak-normalization (see \cite{AndreottiBombieri} Thm.\ 4) that weak-normalization, hence the dual-graph $\Gamma.(Z)$, is functorial with respect to non-constant finite morphisms.
\end{sect}

\begin{sect}\label{para: construction of Z star}
Suppose moreover all points in the normalization $\tilde{Z}$ lying above a singular point $z$ of $Z$ are $k(z)$-rational.
Denote by $\Sigma.(z)$ the star in $\Gamma.(Z)$ of a node corresponding to the singular point $z$ and by $\mathcal{Z}^\ast$ the \'{e}tale homotopy type of the normalization $\tilde{Z}$ glued to $\coprod_zBG_{k(z)} \otimes \Sigma.(z)$ via the points above the singular points $z$ of $Z$.
%Denote by $\Sigma.(Z)$ the disjoint union of the stars $\Sigma.(z)$ in $\Gamma.(Z)$ of nodes corresponding to singular points $z$ and by $\mathcal{Z}^\ast$ the \'{e}tale homotopy type of the normalization $\tilde{Z}$ glued to $BG_K \otimes \Sigma.(Z)$ via the points above singular points of $Z$.
Contracting all the stars $\Sigma.(z)$, we get a canonical factorization of the homotopy type of the normalization $\pi$ over the canonical morphism $\mathcal{Z}^\ast \rightarrow \widehat{\rm Et}(Z/K)$, functorial with respect to non-constant finite morphisms.
\end{sect}

\begin{thm}\label{thm: homotopy type of a curve}
Let $Z$ be a reduced, connected curve over $K$ and assume either $Z$ to be semi-normal or $K$ to be perfect.
%Suppose all the points in the normalization $\pi: \tilde{Z} \rightarrow Z$ lying above singular points are $K$-rational.
Denote by $\pi_i: \tilde{Z}_i \rightarrow Z_i$ the normalizations of the irreducible components and let $\mathcal{R} = \mathcal{R}(Z)$ be the set of indices $i$ s.t.\ $\tilde{Z}_i$ is a rational projective component.
\begin{enumerate}
 \item Suppose all points in the normalization $\tilde{Z}$ lying above a singular point $z$ of $Z$ are $k(z)$-rational.
 Then $\mathcal{Z}^\ast \rightarrow \widehat{\rm Et}(Z/K)$ is a weak equivalence in $\hat{\mathcal{S}}\downarrow BG_K$ (cf.~\ref{para: construction of Z star}).
 \item If $K$ is separably closed, then
 \begin{equation*}
 \widehat{\rm Et}(Z/K) \simeq (\bigvee_i \widehat{\rm Et}(\tilde{Z}_i/K)) \vee \Gamma.(Z)
 \end{equation*}
 holds (functorially in non-constant finite morphisms) in the homotopy category $\mathcal{H}(\hat{\mathcal{S}})$.
 \item Suppose $K$ has characteristic $0$ or cohomological dimension $\leq 1$ and all points in the normalization $\pi: \tilde{Z} \rightarrow Z$ lying above singular points of $Z$ contained in rational projective components are $K$-rational. 
 Then for each rational component $\tilde{Z}_i$ there is a section $s_i$ of $\pi_1(Z/K)$ s.t.\
 \begin{equation*}
  \widehat{\rm Et}(Z/K) \simeq
  (B\pi_1(Z)) \vee_{\coprod_i s_i, \coprod_i x_i} (\coprod_{i\in \mathcal{R}}\widehat{\rm Et}(\tilde{Z}_i/K)) 
 \end{equation*}
 holds in $\mathcal{H}(\hat{\mathcal{S}}\downarrow BG_K)$ for $x_{i}$ an arbitrary $K$-rational point of $\tilde{Z}_{i}$.
\end{enumerate}
\end{thm}

% \begin{rem}\label{rem: strange curves}
% The assumption in Thm.~\ref{thm: homotopy type of a curve} might fail for general curves:
% Let $Z/K$ be a regular curve and $z$ a closed point with $k(z) \neq K$.
% Contracting the Galois-orbit $z\otimes_KK^{\rm s}$ in $Z\otimes_KK^{\rm s}$ gives a curve $Z^{\prime,{\rm s}} /K^{\rm s}$ which descents to a curve $Z^\prime /K$ with normalization $Z$. 
% Under the normalization morphism, $z$ lies over the $K$-rational point of $Z^\prime$ corresponding to the contracted Galois-orbit. 
% \end{rem}

\begin{rem}\label{rem: homotopy type of a curve, non completed case}
If one works with pro-simplicial sets together with Isaksen's model structure (see \cite{Isaksen01}), the proof below for Thm.\ \ref{thm: homotopy type of a curve} (i), for (ii) and for the characteristic $0$ case of (iii) still goes through without profinite completion.
For the remaining case of (iii) an analogue of Lem.~\ref{lem: ho classes of rational points} (i) is needed, e.g.\ an adequate analogue of the Quillen equivalence $\mathcal{H}(\hat{\mathcal{S}}\downarrow B\Gamma) \simeq \mathcal{H}(\hat{\mathcal{S}}_\Gamma)$ in \cite{Quick10} 2.11 and for the descent spectral sequence of loc.\ cit.\ Thm.\ 2.16.
\end{rem}

\begin{rem}\label{rem: homotopy type of a curve, functoriality}
Let $f: Z^\prime \rightarrow Z$ be a finite non-constant morphism of curves as in Thm.\ \ref{thm: homotopy type of a curve} (iii).
Suppose that there is no non-rational projective component $\tilde{Z}_j^\prime$ of $Z^\prime$ lying over a rational projective component $\tilde{Z}_i$ of $Z$.
Then the proof of Thm.\ \ref{thm: homotopy type of a curve} will show that $f$ is compatible with
\begin{equation*}
 (B\pi_1(f)) \vee_{\coprod_i s_i^\prime, \coprod_j x_j^\prime} (\coprod_{j\in \mathcal{R}(Z^\prime)} f\vert_{\tilde{Z}_j^\prime})
\end{equation*}
for compatible choices of the $x_i$ and $s_i$.
If there is a non-rational component $\tilde{Z}_j^\prime$ over a rational projective one $\tilde{Z}_i$, then this is no longer true:
E.g.\ ${\rm H}^2(\tilde{Z}_i) \rightarrow {\rm H}^2(\tilde{Z}_j^\prime)$ is invisible for $B\pi_1(f)$.
\end{rem}

\begin{rem}\label{rem: tubular neighbourhood of a closed point}
For $S$ a $Z$-scheme and $z$ a closed point of $Z$ we just write $S_z^h$ resp.\ $S^\bullet$ for the henselization $S \times_Z {\rm Spec}(\mathcal{O}_{Z,z}^h)$ resp.\ for the punctured scheme $S \times_Z Z \setminus \{ z \}$.
Before giving the proof of Thm.\ \ref{thm: homotopy type of a curve}, let us first recall that the (punctured) \'{e}tale tubular neighbourhood satisfies
\begin{equation*}
 T_{Z,z}^{(\bullet)} \simeq
 %T_{Z_z^{h},z}^{(\bullet)} \simeq
 \widehat{\rm Et}(Z_z^{h,(\bullet)}):
\end{equation*}
We work with Cox's model for the (punctured) tubular neighbourhood in ${\rm Pro}\mathcal{H}(\underline{\rm SSets}_\bullet)$ (see \cite{Cox78.2}).
The non-punctured case is covered by loc.\ cit.\ Thm.\ 2.2.
For the punctured case, use
\begin{equation*}
 T_{Z,z}^{\bullet} = {\rm lim}_{\mathfrak{V}. \in t_{Z,z}^{\rm Ho}} \widehat{\rm Et}(\mathfrak{V}.^{\bullet}) =  {\rm lim}_{(U,u)} {\rm lim}_{\mathfrak{V}. \in t_{U,u}^{\rm Ho}} \widehat{\rm Et}(\mathfrak{V}.^{\bullet}),
\end{equation*}
where $(U,u)$ runs through the strict \'{e}tale neighbourhoods of $z$ and $t_{U,u}^{\rm Ho}$ is the opposite category of the homotopy category of the full subcategory $t_{U,u}$ of degreewise Noetherian and separated simplicial objects $\mathfrak{V}.$ in $U_{\rm  \acute{e}t}$ s.t.\ $u^*\mathfrak{V}.$ is a hypercover of $k(u)$ ($= k(z)$).
We have to compare $\pi_0$, $\pi_1$ and cohomology of the respective homotopy types.
Arguing degreewise and using the descriptions of these homotopy invariants in \cite{Friedlander} Ch.~5 (together with the spectral sequence in loc.~cit.~Prop.~2.4), we may restrict to simplicial schemes of the form $\mathfrak{V}. = {\rm cosk}_m^U\mathfrak{V}.$ in $t_{U,u}$ for $m\gg 0$.
After a possible refinement of $(U,u)$, such an object in $t_{U,u}$ can be refined to a simplicial object $\mathfrak{V}.$ in the finite \'{e}tale site $U_{\rm f\acute{e}t}$.
It follows that the canonical morphisms $\mathfrak{V}_n \rightarrow ({\rm cosk}_{n-1}^U\mathfrak{V}.)_n$ are finite \'{e}tale.
Since $u^*\mathfrak{V}.$ is a hypercover of $k(u)$, these maps are surjective on connected components, i.e., surjective and $\mathfrak{V}.$ is a hypercover of $U$.
In particular, $\mathfrak{V}.^\bullet \rightarrow U^\bullet$ is a hypercover and hence a weak equivalence.
Summing up, we get
\begin{equation*}
 {\rm lim}_{(U,u)} {\rm lim}_{\mathfrak{V}. \in t_{U,u}^{\rm Ho}} \widehat{\rm Et}(\mathfrak{V}.^{\bullet}) \simeq {\rm lim}_{(U,u)} \widehat{\rm Et}(U^{\bullet}) \simeq \widehat{\rm Et}(Z_z^{h,\bullet})
\end{equation*}
and hence the claim follows.
\end{rem}

\noindent \textbf{Proof of Thm.\ \ref{thm: homotopy type of a curve}.}
Weak-normalization is a universal homeomorphism by \cite{AndreottiBombieri} Thm. 4, i.e., we may assume that $Z$ is semi-normal in the perfect case, as well. 
Factor the normalization into finitely many steps
\begin{equation*}
 \tilde{Z} =
 \xymatrix{
  Z^{(n)} \ar[r]^-{\pi^{(n-1)}} &
  Z^{(n-1)} \ar[r]^-{\pi^{(n-2)}} &
  \dots \ar[r]^-{\pi^{(0)}} &
  Z^{(0)}
 } = Z,
\end{equation*}
where in each step we resolve exactly one of the ordinary multiple points.
Say in the $(l+1)^{\rm th}$ step $\pi^{(l)}: Z^{(l+1)} \rightarrow Z^{(l)}$ the multiple point $z^{(l)}$ and let $z_1^{(l)}, \dots, z_{r_l}^{(l)}$ be the $r_l$ distinct points of $Z^{(l+1)}$ over $z^{(l)}$.
We get the homotopy pushout
\begin{equation*}
 \widehat{\rm Et}(Z^{(l+1)}/K) \vee_{z_1^{(l)}, \dots, z_{r_l}^{(l)}} (BG_{k(z^{(l)})} \otimes \Sigma.(z^{(l)}))
\end{equation*}
in $\hat{\mathcal{S}}\downarrow BG_K$ by joining $z_1^{(l)}, \dots, z_{r_l}^{(l)}$ via the ends of the rays of the star $BG_{k(z^{(l)})} \otimes \Sigma.(z^{(l)})$.
By contracting $BG_{k(z^{(l)})} \otimes \Sigma.(z^{(l)})$ to $z^{(l)}$, $\pi^{(l)}$ induces a canonical morphism
\begin{equation*}
 \bar{\pi}^{(l)}:
 \xymatrix{
  \widehat{\rm Et}(Z^{(l+1)}/K) \vee_{z_1^{(l)}, \dots, z_{r_l}^{(l)}} (BG_{k(z^{(l)})} \otimes \Sigma.(z^{(l)})) \ar[r] &
  \widehat{\rm Et}(Z^{(l)}/K)
 }.
\end{equation*}

\noindent \textbf{(i):} It suffices to see that $\bar{\pi}^{(l)}$ is a weak equivalence in $\hat{\mathcal{S}}\downarrow BG_K$.
Clearly, it is enough to deal only with the first step ($l=0$).
To ease the notation, we skip the upper subscripts and write just $Z^\prime$ for $Z^{(1)}$, $z$ for $z^{(0)}$ and $z_1,\dots,z_r$ for $z_1^{(0)}, \dots, z_{r_0}^{(0)}$.
By \cite{Friedlander} Prop.\ 15.6, the canonical morphism $\widehat{\rm Et}(Z^\bullet/K) \vee_{T_{Z,z}^\bullet} T_{Z,z} \rightarrow \widehat{\rm Et}(Z/K)$
is a weak equivalence ($Z^\bullet$ is no longer connected, but the proof goes through without any changes), where the (punctured) tubular neighbourhood $T_{Z,z}^{(\bullet)}$ is weakly equivalent to the \'{e}tale homotopy type of the (punctured) henselization $Z_z^{h,(\bullet)}$ by Rem.\ \ref{rem: tubular neighbourhood of a closed point} (and similar for $Z^\prime$).
By prime avoidance, $Z_z^{h,\bullet}$ is affine of dimension $0$ and hence discrete, where the distinct points are given by the function fields of the components $Z_i$ containing $z$. 
Summing up, we have
\begin{equation*}
 Z^\bullet =
 Z^{\prime,\bullet} =
 Z^\prime \setminus \{z_0,\dots,z_r\},
\end{equation*}
\begin{equation*}
 T_{Z,z} \simeq BG_{k(z)} \simeq
 T_{Z^\prime,z_i},
\end{equation*}
\begin{equation*}
 T_{Z,z}^\bullet \simeq
 \coprod_i T_{Z^\prime,z_i}^\bullet,
\end{equation*}
\begin{equation*}
 \widehat{\rm Et}(Z^\prime/K) \simeq %\setminus \{z_i,\dots,z_r\} \simeq
 (\dots(\widehat{\rm Et}(Z^{\prime,\bullet}/K) \vee_{T_{Z^\prime,z_1}^\bullet} T_{Z^\prime,z_1})\vee_{T_{Z^\prime,z_2}^\bullet} \dots)\vee_{T_{Z^\prime,z_{i-1}}^\bullet} T_{Z^\prime,z_{i-1}},
\end{equation*}
where the last statement is again \cite{Friedlander} Prop.\ 15.6 applied multiple times.
Here we used that all points $z_i$ are $k(z)$-rational by assumption.
We get $\widehat{\rm Et}(Z^\prime/K) \simeq \widehat{\rm Et}(Z^{\prime,\bullet}/K) \vee_{(\coprod_iT_{Z^\prime,z_i}^\bullet)} (\coprod_iT_{Z^\prime,z_i})$,
hence
\begin{align*}
 & \widehat{\rm Et}(Z^\prime/K) \vee_{z_0,\dots z_r} (BG_{k(z)} \otimes\Sigma.(z))
 \\
 \simeq~ & (\widehat{\rm Et}(Z^{\prime,\bullet}/K) \vee_{(\coprod_iT_{Z^\prime,z_i}^\bullet)} (\coprod_iT_{Z^\prime,z_i})) \vee_{(\coprod_iT_{Z^\prime,z_i})} (T_{Z,z}\otimes\Sigma.(z)) 
\\ 
 \simeq~ & \widehat{\rm Et}(Z^\bullet/K) \vee_{T_{Z,z}^\bullet} T_{Z,z},
\end{align*}
i.e.\ $\bar{\pi}^{(0)}$ is a weak equivalence, just as claimed.
If we put all the single steps $\bar{\pi}^{(l)}$ together, we get a weak equivalence in $\hat{\mathcal{S}}\downarrow BG_K$ between the \'{e}tale homotopy types of $Z$ and of its normalization $\tilde{Z}$ with $z_0^{(l)},\dots,z_{r_k}^{(l)}$ joined together as the ends of rays of $\Sigma.(z_l)$, i.e., between $\widehat{\rm Et}(Z/K)$ and $\mathcal{Z}^\ast$.

\noindent \textbf{(ii):}
If $K$ is separably closed, the assumptions of (ii) clearly are satisfied.
Further, by Lem.~\ref{lem: ho classes of rational points} (i), we can move all the $z_j^{(k)}$'s on one component $\tilde{Z}_i$ to a single $K$-point $x_i$ of $\tilde{Z}_i$.
As a result, $\widehat{\rm Et}(Z/K)$ is weakly equivalent to the dual graph $\Gamma.(Z)$ with all the $\widehat{\rm Et}(\tilde{Z}_i/K)$ glued via $x_i$ to the corresponding node of $\Gamma.(Z)$.
Since $\Gamma.(Z)$ is connected, we get the desired equivalence between $\widehat{\rm Et}(Z/K)$ and $(\bigvee_i \widehat{\rm Et}(\tilde{Z}_i/)) \vee \Gamma.(Z)$.

\noindent \textbf{(iii):}
First, note we may assume that $Z$ itself is non-rational.
Since $Z$ is connected, each rational projective component $\tilde{Z}_i$ contains a rational point, i.e., is isomorphic to $\mathbb{P}_K^1$.
\\
Let $Z^\prime \hookrightarrow Z$ be the closed subscheme given by the union of components $Z_i$ for $i\notin \mathcal{R}$ and let $Z^{\prime\prime}\rightarrow Z^\prime$ be the normalization at all singular points of $Z$ contained in rational projective components.
By assumption, all points of $Z^{\prime\prime}$ over such singular points are $K$-rational.
The cohomology resp.\ group-cohomology in degrees $\geq 1$ of one point unions resp.\ free products is the direct sum of the cohomology of the factors.
As a bouquet of $S^1$'s, $\Gamma.(Z^{\prime\prime})$ is a $K(\pi,1)$.
Since the $\tilde{Z}_i$ are $K(\pi,1)$ for every $i\notin \mathcal{R}$, too (see \cite{Stix02} Prop.~A.4.1), $Z^{\prime\prime,{\rm s}}$ is a $K(\pi,1)$.
Hence the same is true for $Z^{\prime\prime}$ itself.
\\
Arguing as in (i), we get that $\widehat{\rm Et}(Z/K)$ is weakly equivalent to the \'{e}tale homotopy type of $Z^{\prime\prime} \amalg \coprod_{i\in\mathcal{R}}\tilde{Z}_i$ glued together via the stars $BG_K\otimes \Sigma.(z)$ corresponding to singular points $z$ in $Z$ contained in rational projective components (cf.~the construction of $\mathcal{Z}^\ast$ in \ref{para: construction of Z star}).
Denote the latter homotopy type by $\mathcal{Z}_{\mathcal{R}}^\ast$, i.e., we have a weak equivalence $\mathcal{Z}_{\mathcal{R}}^\ast \rightarrow \widehat{\rm Et}(Z/K)$.
Let $\mathcal{Z}_{K(\pi,1)}$ be the homotopy type we get from $\mathcal{Z}_{\mathcal{R}}^\ast$ by contracting the \'{e}tale homotopy types of the rational projective components $\tilde{Z}_i$ to $BG_K$ for each $i\in \mathcal{R}$ individually.
It follows from the above discussion, that $\mathcal{Z}_{K(\pi,1)}$ is a $K(\pi,1)$-space, weakly equivalent to $B\pi_1(Z)$.
\\
By Lem.~\ref{lem: ho classes of rational points} we can move all the $z_j^{(l)}$'s on one rational projective component $\tilde{Z}_i$ to a single $K$-point $x_i$ in $\hat{\mathcal{S}}\downarrow BG_K$.
It follows that $\mathcal{Z}_{\mathcal{R}}^\ast$ and $\mathcal{Z}_{K(\pi,1)} \vee_{\coprod_i BG_K} (\coprod_{i\in \mathcal{R}}\widehat{\rm Et}(\tilde{Z}_i/K))$ are weakly equivalent via the $K$-points $x_i$ and gluing-morphisms $BG_K \rightarrow \mathcal{Z}_{K(\pi,1)}$ induced by the contracted components in $\mathcal{Z}_{K(\pi,1)}$.
This translates to our claim, if we let $s_i$ be the sections of $\pi_1(Z/K)$ corresponding to these gluing morphisms under the weak equivalence $\mathcal{Z}_{K(\pi,1)} \simeq B\pi_1(Z)$.
\Endproof

Let us formulate some immediate consequences of Thm.\ \ref{thm: homotopy type of a curve}:

\begin{cor}\label{cor: l-good K(pi,1) property}
Let $Z/K$ be a curve as in Thm.\ \ref{thm: homotopy type of a curve}.
\begin{enumerate}
 \item $Z$ is a $K(\pi,1)$ if and only if all the $\tilde{Z}_i$'s are $K(\pi,1)$, i.e., if none of the $\tilde{Z}_i$'s is a Brauer-Severi curve.
 \item If $K$ is separably closed, then
 \begin{equation*}
  \pi_1(Z) \simeq (\Asterisk_{i\notin \mathcal{R}} \pi_1(\tilde{Z}_i)) \ast F_r.
 \end{equation*}
 In particular, for any prime $\ell$, $\pi_1(Z)$ is $\ell$-good if and only if all the $\pi_1(\tilde{Z}_i)$'s are $\ell$-good.
 \item Suppose $Z/K$ satisfies the assumption of Thm.~\ref{thm: homotopy type of a curve} (iii).
 Then for each constructible local system $\Lambda$ on $Z$, we get an exact triangle of cohomology cochains
 \begin{equation*}
 \xymatrix{
  C^\bullet(Z,\Lambda) \ar[r] &
  C^\bullet(\pi_1(Z),\Lambda) \oplus \bigoplus_{i\in \mathcal{R}} C^\bullet(\tilde{Z}_i,\Lambda) \ar[d]^-{\bigoplus_i(s_i^* + x_i^*)}
 \\
  &\bigoplus_{i\in \mathcal{R}}C^\bullet(G_K,\Lambda) \ar[ul]^-{+1}
 }
\end{equation*}
 \item Moreover, suppose $K$ has cohomological dimension $\leq 1$.
 Then the higher cohomology groups ($q = 2,3$) are given by
 \begin{equation*}
  {\rm H}^q(Z,\Lambda) = {\rm H}^q(\pi_1(Z),\Lambda) \oplus \bigoplus_{i\in\mathcal{R}} {\rm H}^q(\tilde{Z}_i,\Lambda).
 \end{equation*}
\end{enumerate}
\end{cor}

\begin{rem}
Using \cite{Stix06} Ex.\ 5.5, we get the statement corresponding to (iii) for the \'{e}tale fundamental group of $Z$ a semi-stable curve.
\end{rem}

\noindent \textbf{Proof of Cor.\ \ref{cor: l-good K(pi,1) property}.}
To prove (i) it suffices to see that $Z\otimes_K K^{\rm s}$ is a $K(\pi,1)$, i.e., we may assume that $K$ is separably closed.
The cohomology resp.\ group-cohomology in degrees $\geq 1$ of one point unions resp.\ free products is the direct sum of the cohomology of the factors.
As a bouquet of $S^1$'s, $\Gamma.(Z)$ is a $K(\pi,1)$.
Thus, (i) follows from Thm.\ \ref{thm: homotopy type of a curve} (ii).
\\
The first statement of (ii) is again a direct consequence of Thm.\ \ref{thm: homotopy type of a curve} (ii).
For the second statement we argue as in (i) using that pro-$\ell$ completion preserves free products.
\\
For $\mathcal{Z}_{K(\pi,1)} \simeq B\pi_1(Z)$ as in the proof of Thm.\ \ref{thm: homotopy type of a curve} (iii), loc.\ cit.\ gives the homotopy cofibre sequence
\begin{equation*}
 \xymatrix{
  \coprod_{i\in\mathcal{R}}BG_K \ar[rr]^-{(\coprod_i s_i) \amalg (\coprod_i x_i)} &&
  \mathcal{Z}_{K(\pi,1)} \amalg (\coprod_{i\in \mathcal{R}} \widehat{\rm Et}(\tilde{Z}_i/K)) \ar[rr]^-{\rm can} &&
  Z
 }
\end{equation*}
inducing the exact triangle of cohomology cochains, claimed in (iii).
\\
Finally, (iv) is a direct consequence of the long exact sequence of cohomology groups induced by the triangle in (iii). 
\Endproof

\begin{cor}\label{cor: curves have l-good geometric fundamental group}
Let $Z/K$ be a curve as in Thm.\ \ref{thm: homotopy type of a curve} for $K$ separably closed and $\ell \neq {\rm char}(K)$.
Then $\pi_1(Z)$ is $\ell$-good with torsion free abelianization $\pi_1^{\ell, {\rm ab}}(Z)$ of the $\ell$-completion.
\end{cor}

\proof
By Cor.\ \ref{cor: l-good K(pi,1) property} (ii), we may assume that $Z$ is even smooth and projective.
If ${\rm char}(K) = 0$, then $\pi_1(Z)$ is $\ell$-good as finitely generated fundamental group of a Riemann surface.
For the same reason, $\pi_1^{\ell, {\rm ab}}(Z)$ is torsion free.
For ${\rm char}(K)>0$, $\pi_1(Z)$ is $\ell$-good by \cite{Stix02} Prop.\ A.4.1.
If $Z/K$ is proper, there is a smooth lift $\mathfrak{Z} / W(K)$ of $Z$ to the Witt-ring $W(K)$ of $K$ by \cite{SGA1} Exp.\ III Cor.\ 7.4.
Let $Z^\prime$ be a geometric generic fibre $\bar{\eta}^*\mathfrak{Z}$ of $\mathfrak{Z}$.
Then the $\ell$-completed specialization map $\pi_1^\ell(Z^\prime) \rightarrow \pi_1^\ell(Z)$ 
is an isomorphism, which settles the second claim.
In case $X$ is affine, it is a $K(\pi,1)$-space of cohomological dimension $\leq 1$ with $\ell$-good fundamental group.
In particular, $\pi_1^\ell(Z)$ has cohomological dimension $\leq 1$, i.e., is free (see e.g.~\cite{NeukirchSchmidtWingberg} Prop.~3.5.17).
\Endproof

\subsection*{\boldmath{$K(\pi,1)$}-models of \boldmath{$p$}-adic curves}
Let $X$ be a geometrically connected smooth projective curve over the $p$-adic field $k$ admitting a section $s$ of $\pi_1(X/k)$.
In particular, $X$ is a $K(\pi,1)$-space.
We want  to modify $X$ and $s$ until $X$ admits a proper flat model $\eta: X \hookrightarrow \mathfrak{X}$ over $\mathfrak{o}$ that is a $K(\pi,1)$-space, too.
By \cite{SGA1} Exp.~X Thm.~2.1 and proper base change, the \'{e}tale homotopy type of $\mathfrak{X}$ is weakly equivalent to the homotopy type of its reduced special fibre $Y = (\mathfrak{X}\otimes_{\mathfrak{o}}\mathbb{F})_{\rm red}$.
Using Cor.~\ref{cor: l-good K(pi,1) property} (i), $\mathfrak{X}$ is a $K(\pi,1)$ if and only if $Y$ has no rational components.
\\
Recall that a {\it neighbourhood} $(X^\prime,s^\prime)$ of the section $s$ is a finite \'{e}tale covering $f$ together with a compatible section $s^\prime$:
\begin{equation*}
 \xymatrix{
  & \widehat{\rm Et}(X^\prime/k)\phantom{.} \ar[d]^-f
 \\
  BG_k \ar[ur]^-{s^\prime} \ar[r]^-s &
  \widehat{\rm Et}(X/k).
 }
\end{equation*}

\begin{lem}\label{lem: K(pi,1)-model}
Let $X$ be a geometrically connected smooth projective curve over $k$ admitting a section $s$ of $\pi_1(X/k)$.
Then there is a finite extension $k^\prime/k$ and a neighbourhood $(X^\prime,s^\prime)$ of the restricted section $s\otimes_k k^\prime$ of $X\otimes_kk^\prime$ s.t.\ $X^\prime$ admits a (stable) $K(\pi,1)$-model over the normalization $\mathfrak{o}^\prime$ of $\mathfrak{o}$ in $k^\prime/k$.
\end{lem}

\proof
By the Stable Reduction Theorem of Deligne and Mumford (\cite{DeligneMumford69} Cor.\ 2.7) we may assume that $X$ has split stable reduction.
Using \cite{Mochizuki96} Lem.\ 2.9 (see also \cite{PopStix14} Lem.\ 5.3), there is a finite \'{e}tale covering $X^{\prime\prime} \rightarrow X$, s.t. the reduced special fibre of the stable model $\mathfrak{X}^{\prime\prime} / \mathfrak{o}^{\prime\prime}$ of $X^{\prime\prime}$ does not admit any rational components.
Here, $\mathfrak{o}^{\prime\prime}/\mathfrak{o}$ is the normalization inside the field of constants $k^{\prime\prime}$ of $X^{\prime\prime}$, i.e., $X^{\prime\prime}$ is geometrically connected over $k^{\prime\prime}$.
Then $s(G_k) \cap \pi_1(X^{\prime\prime}) \leq s(G_k)$ is open and closed, i.e.\ of the form $s(G_{k^{\prime}})$ for a suitable finite extension $k^{\prime} / k$.
The restriction of $s$ to $G_{k^{\prime}}$ factors through $\pi_1(X^{\prime\prime})$.
Composed with the canonical map to $G_{k^{\prime\prime}}$, we get that $k^{\prime}$ is an extension of $k^{\prime\prime}$.
By construction, $X^{\prime}:= X^{\prime\prime} \otimes_{k^{\prime\prime}} k^{\prime}$ is geometrically connected over $k^{\prime}$ and $s$ restricts to a section of $\pi_1(X^{\prime}/k^{\prime})$, compatible with the restricted section of $\pi_1(X\otimes_kk^{\prime}/k^{\prime})$.
By the Hurwitz formula, the special fibre of $\mathfrak{X}^{\prime\prime} \otimes_{\mathfrak{o}^{\prime\prime}} \mathfrak{o}^\prime$ does not contain any rational components, i.e., $\mathfrak{X}^{\prime\prime} \otimes_{\mathfrak{o}^{\prime\prime}} \mathfrak{o}^\prime$ is a $K(\pi,1)$-model.
Finally, as a base change of a stable model, it is still stable.
\Endproof

\section{Specialized sections and homotopy rational points}\label{sect: specialization}

\subsection*{Specialized homotopy rational points.}
Fix a geometrically connected smooth projective curve $X$ over $k$ of genus $\geq 2$ and a proper flat model $\eta: X \hookrightarrow \mathfrak{X}$ over $\mathfrak{o}$ with reduced special fibre $Y = (\mathfrak{X}\otimes_{\mathfrak{o}}\mathbb{F})_{\rm red}$.  
Say, $\pi_1(X/k)$ admits a section $s$ (i.e., a homotopy rational point $s: BG_k \rightarrow X$), unramified with respect to $\eta$, i.e.,  with trivial ramification map
\begin{equation*}
 {\rm ram}_s:
 \xymatrix{
  I_k \ar@{^(->}[r] \ar@/_1pc/[rrr] &
  G_k \ar[r]^-s &
  \pi_1(X) \ar@{->>}[r]^-{\rm sp} &
  \pi_1(Y)
 }.
\end{equation*}
This is equivalent to $s$ specializing to a section $\bar{s}$ of $\pi_1(Y/\mathbb{F})$, compatible via the specialization maps.

Combining proper base change with \cite{SGA1} Exp.~X Thm.~2.1, we get a specialization map ${\rm sp}: \widehat{\rm Et}(X/k) \rightarrow \widehat{\rm Et}(Y/\mathbb{F})$ of \'{e}tale homotopy types in $\mathcal{H}(\hat{\mathcal{S}}\downarrow BG_{\mathbb{F}})$.
Unfortunately, as soon as $Y$ admits rational components, it is no longer a $K(\pi,1)$-space.
So sections of $\pi_1(Y/\mathbb{F})$ (compatible with section of $\pi_1(X/k)$) do no longer a priori correspond to homotopy rational points $BG_{\mathbb{F}} \rightarrow \widehat{\rm Et}(Y/\mathbb{F})$ over $\mathbb{F}$ (compatible with homotopy rational points of $X$ over $k$ via ${\rm sp}$).
E.g., it is no longer clear if $\bar{s}$ induces a map on cohomology (compatible with $s$).
Thus, we define:

\begin{defn}\label{def: specialized ho rat pt}
Let $X/k$ be a geometrically connected smooth projective curve over a $p$-adic field $k$ of genus $\geq 2$ and $\mathfrak{X}/\mathfrak{o}$ a normal, proper, flat model with reduced special fibre $Y =(\mathfrak{X}\otimes_{\mathfrak{o}}\mathbb{F})_{\rm red}/\mathbb{F}$.
We say that a homotopy rational point $r:BG_k \rightarrow \widehat{\rm Et}(X/k)$ over $k$ \textbf{specializes} to a homotopy rational point $\bar{r}: BG_{\mathbb{F}} \rightarrow \widehat{\rm Et}(Y/\mathbb{F})$ over $\mathbb{F}$, if $r$ and $\bar{r}$ are compatible via the specialization morphism of homotopy types, i.e., if
\begin{equation}\label{eq: def specialized ho rat pt}
 \xymatrix{
  BG_k \ar[r]^-r \ar[d]^-{\rm can.} &
  \widehat{\rm Et}(X/k) \ar[d]^-{\rm sp}
 \\
  BG_{\mathbb{F}} \ar[r]^-{\bar{r}} &
  \widehat{\rm Et}(Y/\mathbb{F})
 }
\end{equation}
commutes in $\mathcal{H}(\hat{\mathcal{S}}\downarrow BG_{\mathbb{F}})$.
We say that $r$ \textbf{specializes in cohomological settings} to $\bar{r}$, if \eqref{eq: def specialized ho rat pt} induces a commutative square of cohomology cochains $C^\bullet(-,\Lambda)$ in the derived category $\mathcal{D}(\underline{\rm Ab})$ for $\Lambda$ any continuous finite $\mathbb{Z}[G_{\mathbb{F}}]$-module.\footnote{Of course, more generally Def.~\ref{def: specialized ho rat pt} makes sense for $X$ a proper $k$-variety.} 
\end{defn}

\begin{rem}\label{rem: def of specialized ho rat pt}
If a homotopy rational point of $X/k$, i.e., a section, specializes to a homotopy rational point of $Y/\mathbb{F}$, the induced section of $\pi_1(Y/\mathbb{F})$ is the unique specialization of the original section. 
In particular, the original section is unramified.
\\
The converse is true if $Y$ is a $K(\pi,1)$-space.
By Lem.~\ref{lem: K(pi,1)-model}, after a base extension $k^\prime/k$ an (unramified) section $s$ has a neighbourhood $(X^\prime,s^\prime)$ admitting a $K(\pi,1)$-model.
However, it is not known if $s^\prime$ or even the base extension $s\otimes_k k^\prime$ is still unramified.
\end{rem}

\subsection*{A candidate for a specialized homotopy rational point.}
Suppose all points of the normalization $\tilde{Y}$ over singular points of $Y$ contained in rational components are $\mathbb{F}$-rational.
This is always true after an unramified base extension $k^\prime/k$.
Note that the section $s\otimes_k k^\prime$ of $X\otimes_kk^\prime$ induced by an unramified section $s$ is still unramified.
%Sections of $\pi_1(Y/\mathbb{F})$ modulo conjugation correspond to homotopy rational points of $B\pi_1(Y)$ over $\mathbb{F}$.
If $Y$ admits rational components, it is no longer a $K(\pi,1)$.
Thus, a priori the homotopy rational points of $\widehat{\rm Et}(Y/\mathbb{F})$ might no longer correspond to sections of $\pi_1(Y/\mathbb{F})$ modulo conjugation.
However, such a canonical correspondence still holds for $Y/\mathbb{F}$ or more generally for any curve $Z/K$ as in Thm.\ \ref{thm: homotopy type of a curve} (iii) and $K$ of cohomological dimension $\leq 1$.
Indeed, this is just Cor.~\ref{cor: sections vs ho rat points}, where a section of $\widehat{\rm Et}(Z/K) \rightarrow B\pi_1(Z)$ is provided by Thm.~\ref{thm: homotopy type of a curve} (iii):

\begin{cor}\label{cor: sections vs ho rational points}
Let $Z/K$ be a curve as in Thm.\ \ref{thm: homotopy type of a curve} (iii) for $K$ of cohomological dimension $\leq 1$.
Then there is a canonical one-to-one correspondence between sections of $\pi_1(Z/K)$ modulo conjugation and homotopy rational points $[BG_K,\widehat{\rm Et}(Z/K)]_{\hat{\mathcal{S}}\downarrow BG_K}$ of $Z$. 
\end{cor}

%\proof
% Since $K$ has cohomological dimension $1$, the injectivity of the canonical map
% \begin{equation}\label{eq: sections vs ho rational points}
%  \xymatrix{
%   [BG_K,\widehat{\rm Et}(Z/K)]_{\hat{\mathcal{S}}\downarrow BG_K} \ar[r]^-{\gamma_*} &
%   [BG_K,B\pi_1(Z)]_{\hat{\mathcal{S}}\downarrow BG_K}
%  }
% \end{equation}
% given by $\gamma: \widehat{\rm Et}(Z/K) \rightarrow B\pi_1(Z)$ follows from Lem.~\ref{lem: injection of hofixed points}.
% The surjectivity of \eqref{eq: sections vs ho rational points} is an immediate consequence of Thm.\ \ref{thm: homotopy type of a curve} (iii):
% From loc.\ cit.\ we get a section $B\pi_1(Z) \rightarrow \widehat{\rm Et}(Z/K)$ of $\gamma$ in $\mathcal{H}(\hat{\mathcal{S}}\downarrow BG_K)$, i.e., \eqref{eq: sections vs ho rational points} is a retraction and hence surjective.
%\Endproof

\begin{rem}\label{rem: sections vs ho rational points}
Thus, for a specialized section $\bar{s}$ there is at least a unique candidate for a specialized homotopy rational point.
In abuse of notation, we will denote it by $\bar{s}$, as well.
Although the compatibility with the specialization maps of homotopy types is unclear for non-$K(\pi,1)$-models, we will show that, as a homotopy rational point, $s$ specializes to $\bar{s}$ in cohomological settings in Thm.~\ref{thm: specialized homotopy rational point}, below.
\end{rem}

\subsection*{Geometric pro-$\boldsymbol{\ell}$ completion and specialization.}
At least, sections always specializes to sections of the geometric pro-$\ell$ completed fundamental group sequence $\pi^{(\ell)}(Y/\mathbb{F})$:

% \begin{sect}\label{para: geometric pro l completion}
% The pro-$\ell$ completion $\pi_1(Y\otimes_{\mathbb{F}}\mathbb{F}^{\rm s}) \rightarrow \pi_1^\ell(Y\otimes_{\mathbb{F}}\mathbb{F}^{\rm s})$ is by definition universal among maps from $\pi_1(Y\otimes_{\mathbb{F}}\mathbb{F}^{\rm s})$ to pro-$\ell$ groups.
% Hence its kernel $\Delta_\ell$  is a characteristic subgroup of $\pi_1(Y\otimes_{\mathbb{F}}\mathbb{F}^{\rm s})$.
% In particular, $\Delta_\ell$ is normal in $\pi_1(Y)$ and we define
% \begin{equation*}
%  \pi_1^{(\ell)}(Y) := \pi_1(Y) / \Delta_\ell
% \end{equation*}
% as the \textbf{geometric pro-\boldmath{$\ell$} completion} of $\pi_1(Y)$.
% It sits in the commutative diagram with exact rows
% \begin{equation}\label{eq: geometric l completion}
%  \xymatrix{
%   \phantom{.}\boldsymbol{1} \ar[r] &
%   \pi_1(Y\otimes_{\mathbb{F}}\mathbb{F}^{\rm s}) \ar[r] \ar@{->>}[d] &
%   \pi_1(Y) \ar[r] \ar@{->>}[d] &
%   G_{\mathbb{F}} \ar[r] \ar@{=}[d] &
%   \boldsymbol{1}\phantom{.}
%  \\
%   \phantom{.}\boldsymbol{1} \ar[r] &
%   \pi_1^\ell(Y\otimes_{\mathbb{F}}\mathbb{F}^{\rm s}) \ar[r] &
%   \pi_1^{(\ell)}(Y) \ar[r] &
%   G_{\mathbb{F}} \ar[r] &
%   \boldsymbol{1}.
%  }
% \end{equation}
% Denote the lower sequence by $\pi_1^{(\ell)}(Y/\mathbb{F})$.
% \end{sect}

\begin{lem}\label{lem: unramified on geometric l completion}
Any (not necessarily unramified) section $s$ of $\pi_1(X/k)$ specializes to a unique section $\bar{s}_\ell$ of $\pi_1^{(\ell)}(Y/\mathbb{F})$, i.e., we have a commutative diagram
\begin{equation*}
 \xymatrix{
  G_k \ar[r]^-s \ar[d]^-{\rm can.} &
  \pi_1(X)\phantom{.} \ar[d]^-{\rm sp}
 \\
  G_{\mathbb{F}} \ar[r]^-{\bar{s}_\ell} &
  \pi_1^{(\ell)}(Y).
 }
\end{equation*}
\end{lem}

\proof
This is a direct consequence of \cite{Stix12} Prop.~91:
The ramification map ${\rm ram}_s$ induces the geometrically pro-$\ell$ completed ramification map ${\rm ram}_s^{(\ell)}: I_k \rightarrow \pi_1^{(\ell)}(Y)$.
By construction, ${\rm ram}_s^{(\ell)}$ factors over $\pi_1^{\ell}(Y\otimes_{\mathbb{F}}\mathbb{F}^{\rm s})$, i.e., ${\rm ram}_s^{(\ell)}$ is trivial by loc.\ cit.\footnote{Note that \cite{Stix12} Prop.\ 91 essentially follows from the Weil conjectures.}.
\Endproof

% \begin{sect}\label{para: quasi geometrically pro l completed homotopy type}
% Suppose all the points in $\tilde{Y}$ lying over singular points of $Y$ contained in rational components are $\mathbb{F}$-rational.
% Then Thm.~\ref{thm: homotopy type of a curve}~(iii) applies and $\widehat{\rm Et}(Y/\mathbb{F})$ is weakly equivalent to the homotopy type $B\pi_1(Y) \vee_{\coprod_i s_i, \coprod_i y_i} (\coprod_{i\in \mathcal{R}} \widehat{\rm Et}(\tilde{Y}_i/\mathbb{F}))$ in $\hat{\mathcal{S}} \downarrow BG_{\mathbb{F}}$ for $s_i$ suitable sections of $\pi_1(Y/\mathbb{F})$ and $y_i$ arbitrary $\mathbb{F}$-rational points in the rational projective components $\tilde{Y}_i$.
% To make use of Lem.\ \ref{lem: unramified on geometric l completion}, we need to modify the $K(\pi,1)$-part:
% The sections $s_i$ induce sections $s_i^{(\ell)}$ of $\pi_1^{(\ell)}(Y/\mathbb{F})$.
% In particular, we get a canonical morphism in $\mathcal{H}(\hat{\mathcal{S}}\downarrow BG_{\mathbb{F}})$
% \begin{equation}\label{eq: geometric l completion on the fet-part}
%  \xymatrix{
%   \widehat{\rm Et}(Y/\mathbb{F}) \ar[r] &
%   B\pi_1^{(\ell)}(Y) \vee_{\coprod_i s_i^{(\ell)}, \coprod_i y_i} (\coprod_{i\in \mathcal{R}} \tilde{Y}_i)
%  } =: \mathcal{Y}^{(\ell)}.
% \end{equation}
% The proof of Cor.\ \ref{cor: sections vs ho rational points} goes through for $\mathcal{Y}^{(\ell)}$, too. 
% \end{sect}

\begin{sect}\label{para: geometrically pro l completed homotopy type}
Suppose all points in $\tilde{Y}$ lying over singular points of $Y$ contained in rational components are $\mathbb{F}$-rational.
Let $\widehat{\rm Et}(Y/\mathbb{F}) \rightarrow \mathcal{Y}^{(\ell)}$ be the geometric pro-$\ell$ completion $\widehat{\rm Et}(Y/\mathbb{F})_{(\ell)}^\wedge$.
By Cor.\ \ref{cor: curves have l-good geometric fundamental group}, $\pi_1(Y\otimes_{\mathbb{F}}\mathbb{F}^{\rm s})$ is $\ell$-good, so $B\pi_1^{(\ell)}(Y) \simeq (B\pi_1(Y))_{(\ell)}^\wedge$ holds by \ref{para: geometric pro l completion of a group}.
In particular, Thm.~\ref{thm: homotopy type of a curve} (iii) gives a section of $\mathcal{Y}^{(\ell)} \rightarrow B\pi_1^{(\ell)}(Y)$.
The proof of Cor.\ \ref{cor: sections vs ho rational points} goes through for $\mathcal{Y}^{(\ell)}$, i.e., there is a canonical one-to-one correspondence between sections of $\pi_1^{(\ell)}(Y/\mathbb{F})$ modulo conjugation and homotopy rational points $[BG_{\mathbb{F}},\mathcal{Y}^{(\ell)}]_{\hat{\mathcal{S}}\downarrow BG_{\mathbb{F}}}$.
\end{sect}

From Lem.\ \ref{lem: unramified on geometric l completion} we get:

\begin{cor}\label{cor: unramified on geometric l completion and ho rat pts}
For any (not necessarily unramified) section $s$ of $\pi_1(X/k)$ there a unique homotopy rational point $\bar{s}_\ell:BG_{\mathbb{F}} \rightarrow \mathcal{Y}^{(\ell)}$ of $\mathcal{Y}^{(\ell)}$ in $\mathcal{H}(\hat{\mathcal{S}}\downarrow BG_{\mathbb{F}})$ inducing the specialized section $\bar{s}_\ell$ of $\pi_1^{(\ell)}(Y/\mathbb{F})$ in Lem.\ \ref{lem: unramified on geometric l completion}.
\end{cor}

%By Cor.\ \ref{cor: curves have l-good geometric fundamental group}, $\pi_1(Y\otimes_{\mathbb{F}}\mathbb{F}^{\rm s})$ is $\ell$-good.
Comparing the respective Hochschild-Serre spectral sequences, we get that $\pi_1(Y) \rightarrow \pi_1^{(\ell)}(Y)$ induces an isomorphism on cohomology with continuous finite $\mathbb{Z}_\ell[G_{\mathbb{F}}]$-module coefficients.
Luckily, the same is true for $\widehat{\rm Et}(Y/\mathbb{F}) \rightarrow \mathcal{Y}^{(\ell)}$:

\begin{lem}\label{lem: cohomology of geometric l completion}
Let $\Lambda$ a continuous finite $\mathbb{Z}_\ell[G_{\mathbb{F}}]$-module.
Then $\widehat{\rm Et}(Y/\mathbb{F}) \rightarrow \mathcal{Y}^{(\ell)}$ induces an isomorphism
\begin{equation*}
 \xymatrix{
  {\rm H}^\bullet(\mathcal{Y}^{(\ell)},\Lambda) \ar[r]^-\sim &
  {\rm H}^\bullet(Y,\Lambda)
 }.
\end{equation*}
\end{lem}

\proof
%$\widehat{\rm Et}(Y/\mathbb{F})\times_{BG_{\mathbb{F}}}EG_{\mathbb{F}}$ resp.\ $\mathcal{Y}^{(\ell)} \times_{BG_{\mathbb{F}}}EG_{\mathbb{F}}$ is weakly equivalent to the one point union of the $K(\pi,1)$-part $B\pi_1(Y\otimes_{\mathbb{F}}\mathbb{F}^{\rm s})$ resp.\ $B\pi_1^\ell(Y\otimes_{\mathbb{F}}\mathbb{F}^{\rm s})$ with the homotopy types of the projective spaces $\tilde{Y}_i \otimes_{\mathbb{F}}\mathbb{F}^{\rm s}$ for $i \in \mathcal{R}$.
$\widehat{\rm Et}(Y/\mathbb{F})\times_{BG_{\mathbb{F}}}EG_{\mathbb{F}} \rightarrow (\widehat{\rm Et}(Y/\mathbb{F})\times_{BG_{\mathbb{F}}}EG_{\mathbb{F}})_\ell^\wedge$ induces an isomorphism on cohomology cochains $C^\bullet(-,\Lambda)$ in $\mathcal{D}^+(\underline{\rm Mod}_{G_{\mathbb{F}}})$ (see \ref{para: geometric pro l completion}).
Thus, it induces an isomorphism between the respective Hochschild-Serre spectral sequences.
\Endproof

% \begin{rem}\label{rem: pro l completion}
% As an alternative to the construction of $\mathcal{Y}^{(\ell)}$, we might also work with the pro-$\ell$ version of Quicks profinite completion applied to $\widehat{\rm Et}(Y/\mathbb{F}) \times_{BG_{\mathbb{F}}}EG_\mathbb{F}$ (cf.~\cite{Quick12} Rem.~3.3). 
% \end{rem}

\begin{rem}\label{rem: quasi pro l completion}
As an alternative to the construction of $\mathcal{Y}^{(\ell)}$, we might also work with the following explicit ``quasi pro-$\ell$ completion'' $\mathcal{Y}_{\rm alt}^{(\ell)}$:
By Thm.~\ref{thm: homotopy type of a curve} (iii), $\widehat{\rm Et}(Y/\mathbb{F})$ is weakly equivalent to $B\pi_1(Y) \vee_{\coprod_i s_i, \coprod_i y_i} (\coprod_{i\in \mathcal{R}} \widehat{\rm Et}(\tilde{Y}_i/\mathbb{F}))$ in $\hat{\mathcal{S}} \downarrow BG_{\mathbb{F}}$ for $s_i$ suitable sections of $\pi_1(Y/\mathbb{F})$ and $y_i$ arbitrary $\mathbb{F}$-rational points in the rational components $\tilde{Y}_i$.
Let $s_i^{(\ell)}$ be the section of $\pi_1^{(\ell)}(Y/\mathbb{F})$ induced by $s_i$.
Then, we get a canonical morphism in $\mathcal{H}(\hat{\mathcal{S}}\downarrow BG_{\mathbb{F}})$
 \begin{equation*}%\label{eq: geometric l completion on the fet-part}
  \xymatrix{
   \widehat{\rm Et}(Y/\mathbb{F}) \ar[r] &
   B\pi_1^{(\ell)}(Y) \vee_{\coprod_i s_i^{(\ell)}, \coprod_i y_i} (\coprod_{i\in \mathcal{R}} \tilde{Y}_i)
  } =: \mathcal{Y}_{\rm alt}^{(\ell)}.
 \end{equation*}
Cor.~\ref{cor: unramified on geometric l completion and ho rat pts} and Lem.~\ref{lem: cohomology of geometric l completion} go through for $\mathcal{Y}_{\rm alt}^{(\ell)}$, too. 
\end{rem}

\subsection*{Application: $\boldsymbol{{\rm cl}_s}$ and specialized sections.}
We want to give a new proof for Prop.~B, i.e., \cite{EsnaultWittenberg09} Cor.\ 3.4, based on the work done so far.\footnote{As Esnault and Wittenberg's original proof, our proof will make essential use of Deligne's theory of weights, too (via Lem.\ \ref{lem: unramified on geometric l completion}).}
Since the cup product ${\rm cl}_s \cup \alpha$ equals $s^*\alpha$ for all classes $\alpha$ in $H^2(X,\mathbb{Z}_\ell(1))$ (see \cite{Stix12} Sect.\ 6.1), Tate-Lichtenbaum duality shows that ${\rm cl}_s$ in ${\rm H}^2(X,\mathbb{Z}_\ell(1))$ lies in the image of the $\ell$-adic Chern class map $\hat{c}_1: {\rm Pic}(X)\otimes\mathbb{Z}_\ell \rightarrow {\rm H}^2(X,\mathbb{Z}_\ell(1))$, if and only if the composition with $s^*$ is the trivial map ${\rm Pic}(X)\otimes\mathbb{Z}_\ell \rightarrow {\rm H}^2(G_k,\mathbb{Z}_\ell(1))$.
Thus, Prop.~B is equivalent to:

\begin{prop}\label{prop: The l-adic chern class mapping is killed by a section}
{\rm (cf.\ \cite{EsnaultWittenberg09} Cor.\ 3.4)}
Let $X$ be a geometrically connected smooth projective curve of genus $\geq 2$ over a $p$-adic field $k$, admitting a section $s$ of $\pi_1(X/k)$.
Then for $[\mathcal{L}] \in {\rm Pic}(X)$ and $\ell \neq p$, $s^*\hat{c}_1[\mathcal{L}] = 0$ in ${\rm H}^2(G_k,\mathbb{Z}_\ell(1))$.
\end{prop}

\proof
Using Lem.~\ref{lem: unramified on geometric l completion}, we get the commutative diagram on cohomology
\begin{equation}\label{eq: fundamental diagram}
 \xymatrix{
  {\rm Pic}(X) \ar[r]^-{\hat{c}_1} &
  {\rm H}^2(X,\mathbb{Z}_\ell(1)) & 
  {\rm H}^2(\pi_1(X),\mathbb{Z}_\ell(1)) \ar[l]_-\cong \ar[r]^-{s^*} &
  {\rm H}^2(G_k,\mathbb{Z}_\ell(1))\phantom{.}
 \\
  {\rm Pic}(\mathfrak{X}) \ar[r]^-{\hat{c}_1} \ar[u]_-{\eta^*} &
  {\rm H}^2(\mathfrak{X},\mathbb{Z}_\ell(1)) \ar[u]_-{\eta^*} \ar[d]^-{\sigma^*} & 
  {\rm H}^2(\pi_1(\mathfrak{X}),\mathbb{Z}_\ell(1)) \ar[l] \ar[u] \ar[d]^-\cong & 
 \\
  &
  {\rm H}^2(Y,\mathbb{Z}_\ell(1)) &
  {\rm H}^2(\pi_1(Y),\mathbb{Z}_\ell(1)) \ar@/_4pc/[uu] \ar[l] &
 \\
  & & 
  {\rm H}^2(\pi_1^{(\ell)}(Y),\mathbb{Z}_\ell(1)) \ar[u]_-\cong \ar[r]^-{\bar{s}_\ell^*} &
  {\rm H}^2(G_{\mathbb{F}},\mathbb{Z}_\ell(1)). \ar[uuu]
}
\end{equation}
Since $\mathbb{F}$ has cohomological dimension $1$, ${\rm H}^2(G_{\mathbb{F}},\mathbb{Z}_\ell(1))$ is trivial.
Thus, Prop.\ \ref{prop: The l-adic chern class mapping is killed by a section} would follow if all horizontal leftward maps were isomorphisms (i.e., if $Y$ is a $K(\pi,1)$-space) and if $[\mathcal{L}]$ lied in the image of $\eta^*$.
To guarantee this, note that we may replace $[\mathcal{L}]$ by a multiple $[\mathcal{L}^{\otimes n}]$ and enlarge the base field $k$ by any finite extension $k^\prime/k$:
Indeed, ${\rm H}^2(G_k,\mathbb{Z}_\ell(1)) = \mathbb{Z}_\ell$ is torsion free, so a non-zero class $s^*\hat{c}_1[\mathcal{L}]$ would stay non-trivial after multiplication by $n$ resp.\ base extension to $k^\prime$.
It is clear that we might replace the pair $(X,s)$ by a neighbourhood $(X^\prime,s^\prime)$.
Thus, by Lem.\ \ref{lem: K(pi,1)-model}, we might indeed assume that $Y$ is a $K(\pi,1)$-space and $\mathfrak{X}$ the stable model of $X$.
As such, $\mathfrak{X}$ has at most rational singularities (this follows from \cite{Lipman69} Thm.\ 27.1, use that we get $\mathfrak{X}$ by contracting rational curves in the minimal regular model with self intersection $-2$).
These are $\mathbb{Q}$-factorial by loc.\ cit.\ \S 17, i.e., for suitable $n \gg 0$, the $n$-th multiple of the closure in $\mathfrak{X}$ of the Weil-divisor corresponding to $[\mathcal{L}]$ is a Cartier-divisor and thus a pre-image for $[\mathcal{L}^{\otimes n}]$ under $\eta^*$.
\Endproof

\begin{rem}\label{rem: specialized horat points whishfull thinking}
Let $\mathfrak{X}$ be a regular (not necessarily $K(\pi,1)$-) model, s.t.\ $Y$ satisfies the assumptions of Thm.~\ref{thm: homotopy type of a curve} (iii).
It will turn out later that the homotopy rational point $s$ ``specializes in cohomological settings'' to the homotopy rational point $\bar{s}_\ell$ of $\mathcal{Y}^{(\ell)}$ given by Cor.\ \ref{cor: unramified on geometric l completion and ho rat pts} (see Cor.\ \ref{cor: specialized homotopy rational point}, below).
Thus, Prop.\ \ref{prop: The l-adic chern class mapping is killed by a section} would follow directly from Lem.\ \ref{lem: cohomology of geometric l completion} and a diagram chase in \eqref{eq: fundamental diagram}.
Unfortunately, we don't know how to prove Cor.\ \ref{cor: specialized homotopy rational point} without using Prop.\ \ref{prop: The l-adic chern class mapping is killed by a section} itself.
\end{rem}

\begin{rem}\label{rem: prop b as specialization result}
Let $f: X \rightarrow \mathbb{P}_k^N$ be a non-constant $k$-morphism, $\mathcal{L} = f^*\mathcal{O}(1)$ and $r = f \circ s$ resp.\ $r_\infty$ the induced homotopy rational point of $\mathbb{P}_k^N /k$ resp.\ $\mathbb{P}_k^\infty /k$ - here we embed $\mathbb{P}_k^N$ into $\mathbb{P}_k^\infty$ via the first $N+1$ coordinate functions.
The homotopy type of $\mathbb{P}_{\mathbb{F}}^\infty$ represents ${\rm H}^2(-,\hat{\mathbb{Z}}^{(p^\prime)}(1))$ in $\mathcal{H}(\hat{\mathcal{S}}\downarrow BG_{\mathbb{F}})$ and the composition of $r_\infty$ with the specialization morphism of homotopy types $\widehat{\rm Et}(\mathbb{P}_k^\infty/k) \rightarrow \widehat{\rm Et}(\mathbb{P}_{\mathbb{F}}^\infty/\mathbb{F})$ corresponds to the class $s^*\hat{c}_1[\mathcal{L}]$ in ${\rm H}^2(G_k,\hat{\mathbb{Z}}^{(p^\prime)}(1))$.
Since ${\rm H}^2(G_{\mathbb{F}},\hat{\mathbb{Z}}^{(p^\prime)}(1))$ is trivial, Prop.~\ref{prop: The l-adic chern class mapping is killed by a section} implies that the homotopy rational point $r_\infty$ specializes to the unique homotopy rational point of $\mathbb{P}_{\mathbb{F}}^\infty$.
The condition that $s^*\hat{c}_1[\mathcal{L}]$ (or more generally $r^*\hat{c}_1[\mathcal{M}]$ for $r$ a homotopy rational point of a Brauer-Severi variety over $k$ and $\mathcal{M}$ a non-trivial line bundle) is trivial even in ${\rm H}^2(G_k,\hat{\mathbb{Z}}(1))$ implies that $r$ is homotopic to a rational point of $\mathbb{P}_k^N$ (see \cite{JSchmidt15} Thm.~4.3).
\end{rem}

\begin{rem}\label{rem: linearized weak section conjecture}
Suppose $s$ trivializes $\hat{c}_1$ even in ${\rm H}^2(G_k,\hat{\mathbb{Z}}(1))$.
By Tate-Lichtenbaum duality, this is equivalent to the algebraicity of the cycle class ${\rm cl}_s$.
In particular, $X$ admits an algebraic cycle of degree $1$ in this case.
In this sense, such a section $s$ satisfies a \textit{linearized form of the $p$-adic weak section conjecture} (which clearly is predicted by the $p$-adic weak section conjecture).
By Roquette-Lichtenbaum, this is equivalent to the triviality of the relative Brauer group ${\rm Br}(X/k)$.
For a general section $s$, Prop.~B implies that ${\rm Br}(X/k)$ is a $p$-group, which was first proved by Stix (see \cite{Stix10} Thm.~15) using different methods and including even the genus $1$ case.
\end{rem}

\section{Specialized sections and cohomology}\label{sect: specialization of ho rational points}

\subsection*{Specialized homotopy rational points in cohomological settings.}
Fix a geometrically connected smooth projective curve $X$ over $k$ of genus $\geq 2$ admitting a section $s$ unramified with respect to the regular proper flat model $\eta: X \hookrightarrow \mathfrak{X}$ over $\mathfrak{o}$ with reduced special fibre $Y = (\mathfrak{X}\otimes_{\mathfrak{o}}\mathbb{F})_{\rm red}$.
Again, we assume that all points of $\tilde{Y}$ lying over singular points of $Y$ contained in rational components are $\mathbb{F}$-rational.
Then $s$ specializes in cohomological settings (the proof will show more generally: in ``sufficiently additive settings'') to the homotopy rational point $\bar{s}$ of $Y/\mathbb{F}$ given by Cor.\ \ref{cor: sections vs ho rational points}:

\begin{thm}\label{thm: specialized homotopy rational point}
Let $X/k$ be a geometrically connected smooth projective curve over a $p$-adic field $k$ of genus $\geq 2$ and $\mathfrak{X}/\mathfrak{o}$ a regular, proper, flat model s.t.\ the reduced special fibre $Y =(\mathfrak{X}\otimes_{\mathfrak{o}}\mathbb{F})_{\rm red}/\mathbb{F}$ satisfies the assumptions of Thm.~\ref{thm: homotopy type of a curve} (iii).
For an unramified section $s$, its specialized section $\bar{s}$ and a constructible $G_{\mathbb{F}}$-module $\Lambda$, \eqref{eq: def specialized ho rat pt} induces a commutative diagram of cohomology cochains in the derived category $\mathcal{D}^+(\underline{\rm Ab})$ resp.\ $\mathcal{D}^+(\underline{\rm Mod}_{G_{\mathbb{F}}})$:
\begin{equation*}
 \xymatrix{
  C^\bullet(X,\Lambda) \ar[r]^-{s^*} &
  C^\bullet(G_k,\Lambda) 
 \\
  C^\bullet(Y,\Lambda) \ar[r]^-{\bar{s}^*} \ar[u]_-{{\rm sp}^*} &
  C^\bullet(G_{\mathbb{F}},\Lambda) \ar[u]
 }
 \xymatrix{
 \ar@{}[d]|*+{\rm ~resp.~}
 \\
 \phantom{a}
 }
 \xymatrix{
  C^\bullet(X \times_{BG_{\mathbb{F}}} EG_{\mathbb{F}},\Lambda) \ar[r]^-{s^*} &
  C^\bullet(BG_k \times_{BG_{\mathbb{F}}} EG_{\mathbb{F}},\Lambda) \phantom{.}
 \\
  C^\bullet(Y \times_{BG_{\mathbb{F}}} EG_{\mathbb{F}},\Lambda) \ar[r]^-{\bar{s}^*} \ar[u]_-{{\rm sp}^*} &
  C^\bullet(EG_{\mathbb{F}},\Lambda). \ar[u]
 }
\end{equation*}
\end{thm}

To prove Thm.\ \ref{thm: specialized homotopy rational point} we try to work out what the obstruction for the commutativity could be:

\begin{sect}\label{para: obstruction for commutativity}
By Thm.\ \ref{thm: homotopy type of a curve} (iii), $\widehat{\rm Et}(Y/\mathbb{F}) \simeq (B\pi_1(Y)) \vee_{\coprod_i s_i, \coprod_i y_i} (\coprod_{i\in \mathcal{R}}\widehat{\rm Et}(\tilde{Y}_i/\mathbb{F}))$ holds in $\mathcal{H}(\hat{\mathcal{S}}\downarrow G_{\mathbb{F}})$ for any $\mathbb{F}$-points $y_i$ of $\tilde{Y}_i \cong \mathbb{P}_{\mathbb{F}}^1$.
Contracting the $K(\pi,1)$-part $B\pi_1(Y)$ in $\hat{\mathcal{S}}\downarrow G_{\mathbb{F}}$ gives a homotopy co-cartesian square
\begin{equation*}
 \xymatrix{
  B\pi_1(Y) \ar[r] \ar[d]^-{\iota} &
  BG_{\mathbb{F}} \ar[d]^-\ast
 \\
  (B\pi_1(Y)) \vee_{\coprod_i s_i, \coprod_i y_i} (\coprod_{i\in \mathcal{R}}\widehat{\rm Et}(\tilde{Y}_i/\mathbb{F})) \ar[r]^-{p} &
  \bigvee_{BG_{\mathbb{F}}}^{i \in \mathcal{R}} (\widehat{\rm Et}(\tilde{Y}_i/\mathbb{F}),y_i)
 }
\end{equation*}
in $\hat{\mathcal{S}}\downarrow BG_{\mathbb{F}}$, i.e., an exact triangle in $\mathcal{D}^+(\underline{\rm Ab})$
\begin{equation}\label{eq: contracting triangle}
 \xymatrix{
  C^\bullet(\bigvee_{BG_{\mathbb{F}}}^{i \in \mathcal{R}} (\widehat{\rm Et}(\tilde{Y}_i/\mathbb{F}),y_i),\Lambda) \ar[rr]^-{p^*\oplus \ast} &&
  C^\bullet(Y,\Lambda) \oplus C^\bullet(G_{\mathbb{F}},\Lambda) \phantom{.} \ar[d]^{\iota^* - {\rm can}}
 \\
  &&
  C^\bullet(\pi_1(Y),\Lambda). \ar[ull]^-{+1}
 }
\end{equation}
The section $s$ together with ${\rm sp}$ and the ``projection'' $p$ induces a morphism
\begin{equation}\label{eq: obstruction class}
 \xymatrix{
  BG_k \ar[r]^-s  \ar@/_1pc/[rrr]&
  \widehat{\rm Et}(X/k) \ar[rr]^-{:= \alpha} &&
  {\rm cosk}_3 \bigvee_{BG_{\mathbb{F}}}^{i \in \mathcal{R}} (\widehat{\rm Et}(\tilde{Y}_i/\mathbb{F}),y_i)
 }
\end{equation}
in $\mathcal{H}(\hat{\mathcal{S}}\downarrow BG_{\mathbb{F}})$.
Assume that it factors as the canonical morphism to $BG_{\mathbb{F}}$ followed by the distinguished point $\ast: BG_{\mathbb{F}} \rightarrow {\rm cosk}_3 \bigvee_{BG_{\mathbb{F}}}^{i \in \mathcal{R}} (\widehat{\rm Et}(\tilde{Y}_i/\mathbb{F},y_i)$.
Truncation in degrees $>3$ gives the commutative diagram
\begin{equation*}
 \xymatrix{
  C^\bullet(\bigvee_{BG_{\mathbb{F}}}^{i \in \mathcal{R}} (\widehat{\rm Et}(\tilde{Y}_i/\mathbb{F}),y_i),\Lambda) \ar[rr] \ar[dr]_-{\ast} &&
  C^\bullet(G_k,\Lambda)
 \\
  & C^\bullet(G_{\mathbb{F}},\Lambda) \ar[ur]_-{\rm can}
 }
\end{equation*}
on cochains:
Indeed, $k$ and $\mathbb{F}$ have cohomological dimension $\leq 2$ and $C^\bullet(\bigvee_{BG_{\mathbb{F}}}^{i \in \mathcal{R}} (\widehat{\rm Et}(\tilde{Y}_i/\mathbb{F}),y_i),\Lambda)$ is the truncation $\tau^{\leq 3} C^\bullet({\rm cosk}_3\bigvee_{BG_{\mathbb{F}}}^{i \in \mathcal{R}} (\widehat{\rm Et}(\tilde{Y}_i/\mathbb{F}),y_i),\Lambda)$, since the rational components geometrically have cohomological dimension $2$.
It follows that
\begin{equation*}
 \xymatrix{
  C^\bullet(Y,\Lambda) \oplus C^\bullet(G_{\mathbb{F}},\Lambda) \ar[rr]^-{({\rm sp}\circ s)^* - {\rm can}} &&
  C^\bullet(G_k,\Lambda)
 }
\end{equation*}
extents to a morphism of the exact triangle (\ref{eq: contracting triangle}) to the trivial triangle
\begin{equation*}
 \xymatrix{
  0 \ar[r] &
  C^\bullet(G_k,\Lambda) \ar[r]^-{\rm id} &
  C^\bullet(G_k,\Lambda) \ar[r]^-{+1} &
  0[1]
 }.
\end{equation*}
In particular, $C^\bullet(Y,\Lambda) \rightarrow C^\bullet(G_k,\Lambda)$ factors through $\iota^*$.
Since $\iota^*$ is a retraction of the canonical morphism $C^\bullet(\pi_1(Y),\Lambda) \rightarrow C^\bullet(Y,\Lambda)$, the resulting morphism $C^\bullet(\pi_1(Y),\Lambda) \rightarrow C^\bullet(G_k,\Lambda)$ is just the morphism given by the original section $s$ and ${\rm sp}$ on fundamental groups.
Combining this with the compatibility of $s$ with the specialized section $\bar{s}$, we get the commutative diagram
\begin{equation*}
 \xymatrix{
  C^\bullet(X,\Lambda) \ar[rr]^-{s^*} &&
  C^\bullet(BG_k,\Lambda) \phantom{,}
 \\
  C^\bullet(Y,\Lambda) \ar[u]_-{{\rm sp}^*} \ar[r]^-{\iota^*} &
  C^\bullet(\pi_1(Y),\Lambda) \ar[r]^-{\bar{s}^*} \ar[ur]^-{({\rm sp}\circ s)^*} &
  C^\bullet(BG_{\mathbb{F}},\Lambda), \ar[u]_-{\rm can}
 }
\end{equation*}
where the ``outer'' commutative square is the square induced by \eqref{eq: def specialized ho rat pt}.
%[Note that that we made these arguments just to avoid using an appropriate derived category for the homology chains of objects in $\hat{\mathcal{S}}_{G_k}$]
Of course, the same arguments work in the $G_{\mathbb{F}}$-equivariant setting, as well.
\end{sect}

Let us sum up \ref{para: obstruction for commutativity}:
To prove Thm.\ \ref{thm: specialized homotopy rational point}, it is enough to show that the composed morphism (\ref{eq: obstruction class}) factors in $\mathcal{H}(\hat{\mathcal{S}}\downarrow BG_{\mathbb{F}})$ as the canonical morphism to $BG_{\mathbb{F}}$ followed by the ``distinguished point'' $\ast: BG_{\mathbb{F}} \rightarrow {\rm cosk}_3 \bigvee_{BG_{\mathbb{F}}}^{i \in \mathcal{R}} (\widehat{\rm Et}(\tilde{Y}_i/\mathbb{F}),y_i)$.
This problem can be translated into the vanishing of a certain cohomology class:

\begin{sect}\label{para: obstruction classes in cohomology}
Under the Quillen-equivalence of $\mathcal{H}(\hat{\mathcal{S}}\downarrow BG_{\mathbb{F}})$ and $\mathcal{H}(\hat{\mathcal{S}}_{G_{\mathbb{F}}})$, ${\rm cosk}_3 \bigvee_{BG_{\mathbb{F}}}^{i \in \mathcal{R}} (\widehat{\rm Et}(\tilde{Y}_i/\mathbb{F}),y_i)$ corresponds to a $K(\bigoplus_{i\in \mathcal{R}}\hat{\mathbb{Z}}^{(p^\prime)}(1),2)$:
$\tilde{Y}_i$ is isomorphic to $\mathbb{P}_{\mathbb{F}}^1$, hence geometrically simply connected and $\mathbb{P}_{\mathbb{F}}^1 \otimes_{\mathbb{F}}\mathbb{F}^{\rm s}$ has cohomological-$p$-dimension $\leq 1$ (in fact $=0$) as projective curve over a separably closed field of characteristic $p$.
Denote by $\alpha_i$ the $i^{\rm th}$-summand in ${\rm H}^2(X,\hat{\mathbb{Z}}(1)^{(p^\prime)})$ of the class corresponding to $\alpha$ in the composition (\ref{eq: obstruction class}).
Since $G_{\mathbb{F}}$ has cohomological dimension $1$, our factorization problem translates to the question whether $s^*\alpha_i$ is trivial in ${\rm H}^2(G_k,\hat{\mathbb{Z}}(1)^{(p^\prime)})$. 
As prime-to-$p$ Tate module of the Brauer group, ${\rm H}^2(G_k,\hat{\mathbb{Z}}(1)^{(p^\prime)})$ is canonically isomorphic to $\hat{\mathbb{Z}}^{(p^\prime)}$ and it remains to show that $s^*\alpha_i$ vanishes in ${\rm H}^2(G_k,\mathbb{Q}_\ell(1))$ for all $i\in\mathcal{R}$ and $\ell \neq p$.
\end{sect}

We want to have a better understanding where the classes $\alpha_i$ come from:

\begin{sect}\label{para: obstruction class induced by the model}
As a morphism in $\mathcal{H}(\hat{\mathcal{S}}\downarrow BG_{\mathbb{F}})$, $\alpha_i$ in ${\rm H}^2(X,\hat{\mathbb{Z}}(1)^{(p^\prime)})$ is given as the $3$-coskeleton of the composition
\begin{equation}\label{eq: obstruction class II}
 \xymatrix{
  \widehat{\rm Et}(X/k) \ar[r]^-{\rm sp} \ar@/_1pc/[rrr] &
  \widehat{\rm Et}(Y/\mathbb{F}) \ar[r]^-p &
  \bigvee_{BG_{\mathbb{F}}}^{j \in \mathcal{R}} (\widehat{\rm Et}(\tilde{Y}_j/\mathbb{F}),y_j) \ar[r] &
  \widehat{\rm Et}(\tilde{Y}_i/\mathbb{F})
 },
\end{equation}
where the last morphism is given by contracting all $\widehat{\rm Et}(\tilde{Y}_j/\mathbb{F})$ for $j\neq i$ to $BG_{\mathbb{F}}$.
Let $\beta_i$ be the class in ${\rm H}^2(Y,\hat{\mathbb{Z}}^{(p^\prime)}(1))$ (resp.\ in ${\rm H}^2(Y,\mathbb{Q}_\ell(1))$) given by the morphism $\widehat{\rm Et}(Y/\mathbb{F}) \rightarrow \widehat{\rm Et}(\tilde{Y}_i/\mathbb{F})$ in between composed with $\widehat{\rm Et}(\tilde{Y}_i/\mathbb{F}) \rightarrow {\rm cosk}_3\widehat{\rm Et}(\tilde{Y}_i/\mathbb{F})$, i.e., $\alpha_i = {\rm sp}^*\beta_i$. 
Since ${\rm cosk}_3 \widehat{\rm Et}(\tilde{Y}_i/\mathbb{F})$ represents ${\rm H}^2(-,\hat{\mathbb{Z}}^{(p^\prime)}(1))$ and ${\rm cosk}_3$ does not change ${\rm H}^2$, $\beta_i$ in ${\rm H}^2(Y,\mathbb{Q}_\ell(1))$ is given as the image of $1$ in ${\rm H}^2(\tilde{Y}_i,\mathbb{Q}_\ell(1)) = \mathbb{Q}_\ell$ under the morphism $\widehat{\rm Et}(Y/\mathbb{F}) \rightarrow \widehat{\rm Et}(\tilde{Y}_i/\mathbb{F})$ in the composition (\ref{eq: obstruction class II}).
\end{sect}

\begin{rem}\label{rem: obstruction class is sharp}
Note that the vanishing of $s^*\alpha_i$ in fact is necessary:
In order for Thm.\ \ref{thm: specialized homotopy rational point} to hold, it is necessary for the class $s^*\alpha_i$ to vanish in ${\rm H}^2(G_k,\mathbb{Q}_\ell(1))$:
By the usual arguments, ${\rm H}^2(G_k,\mathbb{Z}_\ell(1))$ and ${\rm H}^2(X,\mathbb{Z}_\ell(1))$ are given as the limits over the corresponding cohomology groups with coefficients $\boldsymbol{\mu}_{\ell^n}$.
Since $\alpha_i = {\rm sp}^*\beta_i$ hold by \ref{para: obstruction class induced by the model} and $\mathbb{F}$ has cohomological dimension $1$, the triviality of $s^*\alpha_i$ would follow from Thm.\ \ref{thm: specialized homotopy rational point}.
\end{rem}

We want to give a slightly different characterization of the class $\beta_i$ than the one in \ref{para: obstruction class induced by the model}:

\begin{sect}\label{para: obstruction class induced by the model II}
As in the proof of Thm.~\ref{thm: homotopy type of a curve} (iii), let $Y^{\prime} \hookrightarrow Y$ be the union of components $Y_j$ for $j\notin \mathcal{R}$, $Y^{\prime\prime} \rightarrow Y^\prime$ the normalization at all singular points of $Y$ contained in rational components and $\mathcal{Y}_{\mathcal{R}}^\ast$ the homotopy type we get from joining $\widehat{\rm Et}(Y^{\prime\prime}/\mathbb{F})$ to the homotopy types of the remaining rational components $\tilde{Y}_j$ via the stars of the corresponding singularities.
Recall that the projection $\widehat{\rm Et}(Y/\mathbb{F}) \rightarrow \widehat{\rm Et}(\tilde{Y}_i/\mathbb{F})$ is given by $\widehat{\rm Et}(Y/\mathbb{F}) \simeq \mathcal{Y}_{\mathcal{R}}^\ast$, then moving all points in $\tilde{Y}_j$ over singular points on $Y$ to one $\mathbb{F}$-rational point $y_j$ for each $j \in \mathcal{R}$, then contracting the $K(\pi,1)$-part $B\pi_1(Y)$ and finally contracting the other rational components $\widehat{\rm Et}(\tilde{Y}_j/\mathbb{F})$ for $j\neq i$.
\\  
We could as well contract the rational components $\tilde{Y}_j$ for $j\neq i$ and $Y^{\prime\prime}$ in $\mathcal{Y}_{\mathcal{R}}^\ast$ first, then moving all points $z_1,\dots,z_t$ in $\tilde{Y}_i$ over singular points on $Y$ to $y_i$ to get $\widehat{\rm Et}(\tilde{Y}_i/\mathbb{F}) \vee_{y_i} (BG_{\mathbb{F}}\otimes \bar{\Gamma}.(i))$ for a certain contraction $\Gamma.(Y) \twoheadrightarrow \bar{\Gamma}.(i)$ of the dual graph and then contract the graph $\bar{\Gamma}.(i)$ to get the same morphism $\mathcal{Y}_{\mathcal{R}}^\ast \rightarrow \widehat{\rm Et}(\tilde{Y}_i/\mathbb{F})$.
We get the homotopy pushout
\begin{equation*}
 \mathcal{Y}^{(i)} := \widehat{\rm Et}(\tilde{Y}_i/\mathbb{F}) \vee_{z_1,\dots,z_t} (BG_{\mathbb{F}} \otimes \Gamma.(i))
\end{equation*}
after the first step, where $\Gamma.(i)$ is constructed in the same way as $\bar{\Gamma}.(i)$ but without gluing together the rays of stars corresponding to singularities of $Y$ in $Y_i$ and ending in the node corresponding to $Y_i$ and where we glue the remaining ``open'' ends of $BG_{\mathbb{F}}\otimes \Gamma.(i)$ to the corresponding points $z_1,\dots,z_t$ in $\tilde{Y}_i$.
By Lem.~\ref{lem: ho classes of rational points}, $\mathcal{Y}^{(i)}$ and $\widehat{\rm Et}(\tilde{Y}_i/\mathbb{F}) \vee_{y_i} (BG_{\mathbb{F}}\otimes \bar{\Gamma}.(i))$ are weakly equivalent in $\hat{\mathcal{S}}\downarrow BG_{\mathbb{F}}$.
Further, contracting the graph $\bar{\Gamma}.(i)$ induces an isomorphism ${\rm H}^2(\tilde{Y}_i,\mathbb{Q}_\ell(1)) \rightarrow {\rm H}^2(\widehat{\rm Et}(\tilde{Y}_i/\mathbb{F}) \vee_{y_i} (BG_{\mathbb{F}}\otimes \bar{\Gamma}.(i)),\mathbb{Q}_\ell(1))$:
The target is isomorphic to ${\rm H}^2(\tilde{Y}_i,\mathbb{Q}_\ell(1)) \oplus {\rm H}^2(BG_{\mathbb{F}}\otimes \bar{\Gamma}.(i),\mathbb{Q}_\ell(1))$, since the co-cartesian square of $\widehat{\rm Et}(\tilde{Y}_i/\mathbb{F}) \vee_{y_i} (BG_{\mathbb{F}}\otimes \bar{\Gamma}.(i))$ is even homotopy co-cartesian and $\mathbb{F}$ has cohomological dimension $1$. 
For $m$ prime to $p$ we find ${\rm H}^2(BG_{\mathbb{F}}\otimes \bar{\Gamma}.(i),\boldsymbol{\mu}_m) = (\mathbb{F}^\times / (\mathbb{F}^\times)^m )^{\oplus r}$
for $r$ the number of loops in $\Gamma(Y)$ using the Hochschild-Serre spectral sequence.
Hence, ${\rm H}^2(BG_{\mathbb{F}}\otimes \bar{\Gamma}.(i),\mathbb{Z}_\ell(1))$ is torsion.
Summing up, we get the class $\beta_i$ in ${\rm H}^2(Y,\mathbb{Q}_\ell(1))$ as the image of $1$ under the map
\begin{equation}\label{eq: obstruction class III}
 \mathbb{Q}_\ell =
 \xymatrix{
  {\rm H}^2(\mathcal{Y}^{(i)},\mathbb{Q}_\ell(1)) \ar[r] &
  {\rm H}^2(Y,\mathbb{Q}_\ell(1))
 }.
\end{equation}
\end{sect}

Having made all these preparations, we get Thm.\ \ref{thm: specialized homotopy rational point} as a consequence of Prop.\ \ref{prop: The l-adic chern class mapping is killed by a section}:

\phantom{leere Zeile}

\noindent \textbf{Proof of Thm.\ \ref{thm: specialized homotopy rational point}.}
For each rational component $\tilde{Y}_i$ of $\tilde{Y}$, choose a closed point $x_i$ of $X$ whose specialization lies in the component $Y_i$ outside the singular locus of $Y$ and $\neq y_i$.
We claim that $\beta$ generates the same $\mathbb{Q}_\ell$-subspace in ${\rm H}^2(\mathfrak{X},\mathbb{Q}_\ell(1))$ than the Chern class $\hat{c}_1[\mathcal{O}_{\mathfrak{X}}(D_i)]$ for $D_i = D$ the divisor we get as the closure of $x_i$ in $\mathfrak{X}$:
\\
By Lem.\ \ref{lem: computing continuous l-adic cohomology}, below, ${\rm H}_D^2(\mathfrak{X},\mathbb{Q}_\ell(1))$ can be computed as continuous cohomology of the $BG_{\mathbb{F}}$-space $\widehat{\rm Et}(\mathfrak{X}/\mathfrak{o})\vee_{\widehat{\rm Et}(\mathfrak{X} - D/\mathfrak{o})}BG_{\mathbb{F}}$ (note that $\widehat{\rm Et}(\mathfrak{o},\bar{\sigma}) = BG_{\mathbb{F}}$ since each hypercovering in the \'{e}tale cite of $\mathfrak{o}$ can be refined to a hypercovering in the finite \'{e}tale site).
Let $Y^\bullet$ be the punctured curve $Y - {\rm sp}(x_i)$ and $\mathcal{Y}_{\mathcal{R}}^{\bullet,\ast} \simeq \widehat{\rm Et}(Y^\bullet/\mathbb{F})$ the analogue construction to $\mathcal{Y}_{\mathcal{R}}^\ast$ with $\tilde{Y}_i$ replaced by the punctured component $\tilde{Y}_i - {\rm sp}(x_i)$.\footnote{Note that $\mathcal{Y}_{\mathcal{R}}^{\bullet,\ast} \neq \mathcal{Z}_{\mathcal{R}(Z)}^\ast$ for $Z = Y^\bullet$ in the proof of Thm.~\ref{thm: homotopy type of a curve} (iii): $Y_i - {\rm sp}(x_i)$ is not projective so it would be a component of $Z^\prime$.}
Then $\widehat{\rm Et}(Y/\mathbb{F})\vee_{\widehat{\rm Et}(Y^\bullet/\mathbb{F})}BG_{\mathbb{F}}$ is isomorphic to $\mathcal{Y}_\mathcal{R}^\ast\vee_{\mathcal{Y}_{\mathcal{R}}^{\bullet,\ast}}BG_{\mathbb{F}}$ by the proof of Thm.\ \ref{thm: homotopy type of a curve} (iii).
Similar to $\mathcal{Y}^{(i)}$, we get
\begin{equation*}
 \mathcal{Y}^{(i),\bullet} := \widehat{\rm Et}(\tilde{Y}_i - {\rm sp}(x_i)/\mathbb{F}) \vee_{z_1,\dots,z_t} (BG_{\mathbb{F}} \otimes \Gamma.(i))
\end{equation*}
after contracting all rational components $\tilde{Y}_j$ for $j\neq i$ and $Y^{\prime\prime}$ in $\mathcal{Y}_{\mathcal{R}}^{\bullet,\ast}$.
Summing up, we get a commutative square 
\begin{equation*}
 \xymatrix{
  \widehat{\rm Et}(\mathfrak{X} - D/\mathfrak{o}) \ar@{^(->}[r] &
  \widehat{\rm Et}(\mathfrak{X}/\mathfrak{o}) \ar[r] &
  \widehat{\rm Et}(\mathfrak{X}/\mathfrak{o}) \vee_{\widehat{\rm Et}(\mathfrak{X} - D/\mathfrak{o})} BG_{\mathbb{F}}
 \\
  \mathcal{Y}_{\mathcal{R}}^{\bullet,\ast} \ar@{^(->}[r] \ar[u] \ar[d] &
  \mathcal{Y}_{\mathcal{R}}^{\ast} \ar[r] \ar[u]^-\sim \ar[d] &
  \mathcal{Y}_{\mathcal{R}}^{\ast} \vee_{\mathcal{Y}_{\mathcal{R}}^{\bullet,\ast}} BG_{\mathbb{F}} \ar[u] \ar[d]_-\sim
 \\
  \mathcal{Y}^{(i),\bullet} \ar@{^(->}[r] &
  \mathcal{Y}^{(i)} \ar[r] &
  \mathcal{Y}^{(i)} \vee_{\mathcal{Y}^{(i),\bullet}} BG_{\mathbb{F}}
 }
\end{equation*}
in $\mathcal{H}(\hat{\mathcal{S}}\downarrow BG_{\mathbb{F}})$.
Observe that $Y^{\prime\prime}$ and all the rational components $\tilde{Y}_j$ for $j\neq i$ we contracted in $\mathcal{Y}^{(i)}$ resp.\ $\mathcal{Y}^{(i),\bullet}$ factor through $\mathcal{Y}_{\mathcal{R}}^{\bullet,\ast}$, so the lower right vertical morphism is a weak equivalence.
It follows that the canonical morphism $\widehat{\rm Et}(\mathfrak{X}/\mathfrak{o}) \rightarrow \widehat{\rm Et}(\mathfrak{X}/\mathfrak{o}) \vee_{\widehat{\rm Et}(\mathfrak{X} - D/\mathfrak{o})} BG_{\mathbb{F}}$ factors through $\mathcal{Y}_{\mathcal{R}}^\ast \rightarrow \mathcal{Y}^{(i)}$, i.e., the (non-trivial) cycle class map ${\rm H}^2_D(\mathfrak{X},\mathbb{Q}_\ell(1)) \rightarrow {\rm H}^2(\mathfrak{X},\mathbb{Q}_\ell(1))$ factors through (\ref{eq: obstruction class III}), whose image is generated by $\beta_i$.
\\
Now $\eta^*\hat{c}_1[\mathcal{O}_{\mathfrak{X}}(D_i)] = \hat{c}_1[\mathcal{O}_X(x_i)]$ lies in the kernel of $s^*$ by Prop.\ \ref{prop: The l-adic chern class mapping is killed by a section}.
Thus, the same is true for $\alpha_i$ and Thm.\ \ref{thm: specialized homotopy rational point} follows from \ref{para: obstruction classes in cohomology} and \ref{para: obstruction for commutativity}.
\Endproof

\begin{lem}\label{lem: computing continuous l-adic cohomology}
Let $\mathfrak{X}$ be as in Thm.~\ref{thm: specialized homotopy rational point}, $\ell \neq p$ a prime and $D$ a regular prime-divisor of $\mathfrak{X}$.
Then the continuous cohomology groups ${\rm H}^2(\widehat{\rm Et}(\mathfrak{X}/\mathfrak{o})\vee_{\widehat{\rm Et}(\mathfrak{X} - D/\mathfrak{o})}BG_{\mathbb{F}},\mathbb{Q}_\ell(1))$ and ${\rm H}_D^2(\mathfrak{X},\mathbb{Q}_\ell(1))$ are canonically isomorphic.
\end{lem}

\proof
For $m$ prime to $p$ take the (non-canonically given) morphism of exact triangles
\begin{equation*}
 \xymatrix{
  C^\bullet(\mathfrak{X},\mathfrak{X}-D,\boldsymbol{\mu}_m) \ar[r] \ar[d]&
  C^\bullet(\mathfrak{X},\boldsymbol{\mu}_m) \ar[r] \ar@{^(->}[d]^-{\rm can}&
  C^\bullet(\mathfrak{X}-D,\boldsymbol{\mu}_m) \ar[r]^-{+1} \ar@{=}[d]&
 \\
  C^\bullet(\widehat{\rm Et}(\mathfrak{X}/\mathfrak{o})\vee_{\widehat{\rm Et}(\mathfrak{X} - D/\mathfrak{o})}BG_{\mathbb{F}},\boldsymbol{\mu}_m) \ar[r] &
  C^\bullet(\mathfrak{X},\mu_m) \oplus C^\bullet(G_{\mathbb{F}},\boldsymbol{\mu}_m) \ar[r] &
  C^\bullet(\mathfrak{X}-D,\boldsymbol{\mu}_m) \ar[r]^-{+1} &
 }
\end{equation*}
where the upper triangle computes ${\rm H}_D^\bullet(\mathfrak{X},\boldsymbol{\mu}_m)$ by \cite{Friedlander} Cor.\ 14.5 and Prop.\ 14.6 and the lower triangle is given by the homotopy co-cartesian square
\begin{equation*}
 \xymatrix{
  \widehat{\rm Et}(\mathfrak{X}-D/\mathfrak{o}) \ar[r] \ar[d] &
  BG_{\mathbb{F}} \ar[d]
 \\
  \widehat{\rm Et}(\mathfrak{X}/\mathfrak{o}) \ar[r] &
  \widehat{\rm Et}(\mathfrak{X}/\mathfrak{o})\vee_{\widehat{\rm Et}(\mathfrak{X} - D/\mathfrak{o})}BG_{\mathbb{F}}
 }
\end{equation*}
in $\hat{\mathcal{S}}\downarrow BG_{\mathbb{F}}$.
Note that the induced maps ${\rm H}_D^q(\mathfrak{X},\boldsymbol{\mu}_m) \rightarrow {\rm H}^q(\widehat{\rm Et}(\mathfrak{X}/\mathfrak{o})\vee_{\widehat{\rm Et}(\mathfrak{X} - D/\mathfrak{o})}BG_{\mathbb{F}},\boldsymbol{\mu}_m)$ are uniquely determined up to canonical isomorphisms.
The relevant cohomology groups (i.e.\ ${\rm H}^1$ for $\mathfrak{X}$, $\mathfrak{X} - D$, $G_{\mathbb{F}}$ and ${\rm H}^2$ for $\mathfrak{X}$, $G_{\mathbb{F}}$) are all finite by proper base change and \cite{SGA4.5} Th.\ finitude together with the finiteness for $G_{\mathbb{F}}$-cohomology of finite prime-to-$p$ modules, resp.\ for ${\rm H}^1(\mathfrak{X}-D,\boldsymbol{\mu}_m)$ by \cite{SGA4.5} Cycle Prop.\ 2.1.4.
It follows that the resulting short exact sequences computing ${\rm H}_D^2(\mathfrak{X},\boldsymbol{\mu}_m)$ resp.\ ${\rm H}^2(\widehat{\rm Et}(\mathfrak{X}/\mathfrak{o})\vee_{\widehat{\rm Et}(\mathfrak{X} - D/\mathfrak{o})}BG_{\mathbb{F}},\boldsymbol{\mu}_m)$ stay exact after taking the limit over all $m=\ell^n$ and compute the respective continuous cohomology groups (by \cite{Jannsen88} (3.1)).
Everything stays exact after $(-)\otimes_{\mathbb{Z}_\ell} \mathbb{Q}_\ell$ and the resulting canonical map ${\rm H}_D^2(\mathfrak{X},\mathbb{Q}_\ell(1)) \rightarrow {\rm H}^2(\widehat{\rm Et}(\mathfrak{X}/\mathfrak{o})\vee_{\widehat{\rm Et}(\mathfrak{X} - D/\mathfrak{o})}BG_{\mathbb{F}},\mathbb{Q}_\ell(1))$ is an isomorphism:
${\rm H}^2(G_{\mathbb{F}}, \mathbb{Z}_\ell(1))$ is trivial and ${\rm H}^1(G_{\mathbb{F}}, \mathbb{Z}_\ell(1))$ is a quotient of the multiplicative group $\mathbb{F}^\times$ and hence torsion, since ${\rm H}^1(G_{\mathbb{F}}, \boldsymbol{\mu}_m) = \mathbb{F}^\times/(\mathbb{F}^\times)^m$.
\Endproof

\begin{rem}\label{rem: factorization result without unramified hypothesis}
Looking back onto the proof resp.\ the reductions leading to the proof of Thm.\ \ref{thm: specialized homotopy rational point}, we see that essentially we just proved the factorization
\begin{equation*}
 \phantom{C^\bullet(\pi_1(X),\Lambda) =}
 \xymatrix{
  C^\bullet(X,\Lambda) \ar[r]^-{s^*} &
  C^\bullet(BG_k,\Lambda)
 \\
  C^\bullet(Y,\Lambda) \ar[u]_-{{\rm sp}^*} \ar[r]^-{\iota^*} &
  C^\bullet(\pi_1(Y),\Lambda) \ar[u]_-{({\rm sp}\circ s)^*}
 }
\end{equation*}
for any (possibly ramified) section $s$ of $\pi_1(X/k)$, where $\iota$ was the homotopy splitting of the canonical morphism $\widehat{\rm Et}(Y/\mathbb{F}) \rightarrow B\pi_1(Y)$ given by Thm.\ \ref{thm: homotopy type of a curve} (iii).
\end{rem}

Combining Cor.~\ref{cor: unramified on geometric l completion and ho rat pts} and Lem.~\ref{lem: cohomology of geometric l completion} with the observation in Rem.\ \ref{rem: factorization result without unramified hypothesis} implies:

\begin{cor}\label{cor: specialized homotopy rational point}
Let $X/k$ be a geometrically connected smooth projective curve over a $p$-adic field $k$ of genus $\geq 2$, $\mathfrak{X}/\mathfrak{o}$ a regular, proper, flat model s.t.\ the reduced special fibre $Y/\mathbb{F}$ satisfies the assumptions of Thm.\ \ref{thm: homotopy type of a curve} (iii).
Then for any section $s$ of $\pi_1(X/k)$
there is a unique homotopy rational point $\bar{s}_\ell$ in $[BG_{\mathbb{F}},\mathcal{Y}^{(\ell)}]_{\hat{\mathcal{S}}\downarrow BG_{\mathbb{F}}}$ inducing commutative diagrams 
\begin{equation*}
 \xymatrix{
  G_k \ar[r]^-s \ar[d] &
  \pi_1(X) \ar[d]^-{\rm sp}
 \\
  G_{\mathbb{F}} \ar[r]^-{\bar{s}_\ell} &
  \pi_1^{(\ell)}(Y)
 }
 \xymatrix{
 \ar@{}[d]|*+{\rm ~and~}
 \\
 \phantom{a}
 }
 \xymatrix{
  C^\bullet(X,\Lambda) \ar[r]^-{s^*} &
  C^\bullet(G_k,\Lambda)
 \\
  C^\bullet(Y,\Lambda) \ar[r]^-{\bar{s}_\ell^*} \ar[u]_-{{\rm sp}^*} &
  C^\bullet(G_{\mathbb{F}},\Lambda) \ar[u]
 }
 \xymatrix{
 \ar@{}[d]|*+{\rm ~resp.~}
 \\
 \phantom{a}
 }
 \xymatrix{
  C^\bullet(X\otimes_kk^{\rm nr},\Lambda) \ar[r]^-{s^*} &
  C^\bullet({\rm I}_k,\Lambda)
 \\
  C^\bullet(Y^{\rm s},\Lambda) \ar[r]^-{\bar{s}_\ell^*} \ar[u]_-{{\rm sp}^*} &
  \Lambda \ar[u]
 }
\end{equation*}
in the category of profinite groups resp.\ $\mathcal{D}^+(\underline{\rm Ab})$ resp.\ $\mathcal{D}^+(\underline{\rm Mod}_{G_{\mathbb{F}}})$ for $\Lambda$ any continuous finite $\mathbb{Z}_\ell[G_{\mathbb{F}}]$-module. \Endproof
\end{cor}

\subsection*{Application: A canonical lift of $\boldsymbol{{\rm cl}_s}$.}
In \cite{EsnaultWittenberg09} Rem.\ 3.7 (iii) Esnault and Wittenberg raised the question whether the $\ell$-adic cycle class ${\rm cl}_s$ of a section $s$ of $\pi_1(X/k)$ admits a canonical lift to ${\rm H}^2(\mathfrak{X},\mathbb{Z}_\ell(1))$.
This is predicted by the $p$-adic section conjecture:
If $s$ comes from a rational point $x$, ${\rm cl}_s$ is given by the Chern class of the corresponding divisor $D = x$ and the Chern class of its closure is a canonical lift.
As an application of Cor.\ \ref{cor: specialized homotopy rational point}, we will sketch the construction of a canonical lift of ${\rm cl}_s$ to ${\rm H}^2(\mathfrak{X},\mathbb{Z}_\ell(1))$ for $\mathfrak{X}/\mathfrak{o}$ a regular, proper, flat model s.t.\ the reduced special fibre satisfies the assumptions of Thm.\ \ref{thm: homotopy type of a curve} (iii).
First, let us shortly recall the definition of the cycle class ${\rm cl}_s$:

\begin{rem}\label{rem: cycle class}
Let $X_s \rightarrow X$ be the universal neighbourhood of the section $s$ and $\Delta_X$ the diagonal in $X\times_kX$.
Then ${\rm H}^2(X,\boldsymbol{\mu}_{\ell^n})$ is isomorphic to ${\rm H}^2(X\times_kX_s,\boldsymbol{\mu}_{\ell^n})$ and ${\rm cl}_s$ corresponds to the pullback of $\hat{c}_1[\Delta_X]$ along $X\times_kX_s \rightarrow X\times_kX$.
Since $X$ is a $K(\pi,1)$, $X_s \rightarrow X$ is weakly equivalent to $s$ as a homotopy rational point.
In particular, we get the pullback map as $G_k$-hypercohomology of
\begin{equation*}
 {\rm id}_X^* \otimes^{\mathbb{L}} s^*:
 \xymatrix{
  C^\bullet(X^{\rm s}, \boldsymbol{\mu}_{\ell^n}) \otimes^{\mathbb{L}} C^\bullet(X^{\rm s}, \mathbb{Z} / \ell^n) \ar[r] &
  C^\bullet(X^{\rm s}, \boldsymbol{\mu}_{\ell^n}) \otimes^{\mathbb{L}} C^\bullet(EG_k, \mathbb{Z} / \ell^n)
 }
 \simeq C^\bullet(X^{\rm s}, \boldsymbol{\mu}_{\ell^n}).
\end{equation*}
\end{rem}

Let us first provide a lift of the pullback map ${\rm id}_X^* \otimes^{\mathbb{L}} s^*$, or more generally, of the analogue pullback map ${\rm id}_U^* \otimes^{\mathbb{L}} s^*: C^\bullet(U\times_kX,\boldsymbol{\mu}_{\ell^n}) \rightarrow C^\bullet(U,\boldsymbol{\mu}_{\ell^n})$ for an open subscheme $U \hookrightarrow X$:

\begin{lem}\label{lem: canonical lift of the pullback map}
Let $X/k$ be a geometrically connected smooth projective curve over a $p$-adic field $k$ of genus $\geq 2$, $\mathfrak{X}/\mathfrak{o}$ a regular, proper, flat model s.t.\ the reduced special fibre $Y=(\mathfrak{X}\otimes_{\mathfrak{o}}\mathbb{F})_{\rm red}/\mathbb{F}$ satisfies the assumptions of Thm.\ \ref{thm: homotopy type of a curve} (iii).
Let $\mathfrak{U}\hookrightarrow \mathfrak{X}$ be an open subscheme with generic fibre $U\hookrightarrow X$.
Then for any $\ell \neq p$ and any section $s$ of $\pi_1(X/k)$, there is a canonical lift ${\rm id}_{\mathfrak{U}}^* \otimes^{\mathbb{L}} \bar{s}_\ell^*$ of ${\rm id}_U^* \otimes^{\mathbb{L}} s^*$ to $\mathfrak{U}\times_{\mathfrak{o}}\mathfrak{X}$:
\begin{equation*}
 \xymatrix{
  C^\bullet(U\times_kX,\boldsymbol{\mu}_{\ell^n}) \ar[rr]^-{{\rm id}_U^* \otimes^{\mathbb{L}} s^*} &&
  C^\bullet(U,\boldsymbol{\mu}_{\ell^n}) \phantom{.}
 \\
  C^\bullet(\mathfrak{U}\times_{\mathfrak{o}}\mathfrak{X},\boldsymbol{\mu}_{\ell^n})\ar[rr]^-{{\rm id}_{\mathfrak{U}}^* \otimes^{\mathbb{L}} \bar{s}_\ell^*} \ar[u]_-{\eta^\ast} &&
  C^\bullet(\mathfrak{U},\boldsymbol{\mu}_{\ell^n}) \ar[u]_-{\eta^\ast}.
 }
\end{equation*}
\end{lem}

\proof
From $\bar{s}_\ell^*: C^\bullet(\mathfrak{X}^{\rm nr},\mathbb{Z} / \ell^n)\simeq C^\bullet(Y^{\rm s},\mathbb{Z} / \ell^n) \rightarrow C^\bullet(EG_{\mathbb{F}},\mathbb{Z} / \ell^n)$ in $\mathcal{D}^+(\underline{\rm Mod}_{G_{\mathbb{F}}})$ as in Cor.\ \ref{cor: specialized homotopy rational point}, we get
\begin{equation*}
 {\rm id}_{\mathfrak{U}}^* \otimes^{\mathbb{L}} \bar{s}_\ell^*:
 \xymatrix{
  C^\bullet(\mathfrak{U}^{\rm nr}, \boldsymbol{\mu}_{\ell^n}) \otimes^{\mathbb{L}} C^\bullet(\mathfrak{X}^{\rm nr}, \mathbb{Z} / \ell^n) \ar[r] &
  C^\bullet(\mathfrak{U}^{\rm nr}, \boldsymbol{\mu}_{\ell^n}) \otimes^{\mathbb{L}} C^\bullet(EG_{\mathbb{F}}, \mathbb{Z} / \ell^n)
 }
 \simeq C^\bullet(\mathfrak{U}^{\rm nr}, \boldsymbol{\mu}_{\ell^n}).
\end{equation*}
Taking $G_{\mathbb{F}}$-hypercohomology, we get a canonical pullback map $C^\bullet(\mathfrak{U}\times_{\mathfrak{o}}\mathfrak{X},\boldsymbol{\mu}_{\ell^n}) \rightarrow C^\bullet(\mathfrak{U},\boldsymbol{\mu}_{\ell^n})$.
Consider the commutative diagram
\begin{equation*}
 \xymatrix{
  C^\bullet(\mathfrak{U}\times_{\mathfrak{o}}\mathfrak{X},\boldsymbol{\mu}_{\ell^n}) \ar[r]^-{{\rm id}_{\mathfrak{U}}^* \otimes^{\mathbb{L}} \bar{s}_\ell^*} \ar[d] &
  C^\bullet(\mathfrak{U},\boldsymbol{\mu}_{\ell^n}) \phantom{:} \ar[d]
 \\
  \mathbb{R}\Gamma(G_{\mathbb{F}}, C^\bullet(U^{\rm nr},\boldsymbol{\mu}_{\ell^n}) \otimes^{\mathbb{L}} C^\bullet(X^{\rm nr},\mathbb{Z} / \ell^n)) \ar[r] \ar[d] &
  \mathbb{R}\Gamma(G_{\mathbb{F}}, C^\bullet(U^{\rm nr},\boldsymbol{\mu}_{\ell^n}) \otimes^{\mathbb{L}} C^\bullet({\rm I}_k,\mathbb{Z} / \ell^n)) \phantom{:} \ar[d]
 \\
  \mathbb{R}\Gamma(G_{k}, C^\bullet(U^{\rm nr},\boldsymbol{\mu}_{\ell^n}) \otimes^{\mathbb{L}} C^\bullet(X^{\rm nr},\mathbb{Z} / \ell^n)) \ar[r] \ar[d] &
  \mathbb{R}\Gamma(G_{k}, C^\bullet(U^{\rm nr},\boldsymbol{\mu}_{\ell^n}) \otimes^{\mathbb{L}} C^\bullet({\rm I}_k,\mathbb{Z} / \ell^n)) \phantom{:} \ar[d]
 \\
  C^\bullet(U\times_kX,\boldsymbol{\mu}_{\ell^n}) \ar[r]^-{{\rm id}_U^* \otimes^{\mathbb{L}} s^*} &
  C^\bullet(U,\boldsymbol{\mu}_{\ell^n}):
 }
\end{equation*}
The upper square is induced on $G_{\mathbb{F}}$-hypercohomology by the $\mathcal{D}^+(\underline{\rm Mod}_{G_{\mathbb{F}}})$-square in Cor.\ \ref{cor: specialized homotopy rational point} tensored with $\eta^\ast:C^\bullet(\mathfrak{U}^{\rm nr},\boldsymbol{\mu}_{\ell^n}) \rightarrow C^\bullet(U^{\rm nr},\boldsymbol{\mu}_{\ell^n})$.
The middle square is induced by $G_k \twoheadrightarrow G_{\mathbb{F}}$ (use that an object is flat in $\underline{\rm Mod}_{\mathbb{Z}/\ell^n\mathbb{Z}[\Gamma]}$ if and only if it is flat as a $\mathbb{Z} / \ell^n\mathbb{Z}$-module).
Finally, the lower square is induced on $G_k$-hypercohomology by the $G_k$-equivariant coverings $U^{\rm s} \rightarrow U^{\rm nr}$ resp.\ $EG_k \rightarrow BG_k \times_{BG_{\mathbb{F}}} EG_{\mathbb{F}} \simeq B{\rm I}_k$.
It follows that the compositions of the left resp.\ right vertical maps are the maps on cohomology induced by the respective immersions of the generic fibres.
\Endproof

To get a canonical lift of the $\ell$-adic cycle class map ${\rm cl}_s$, we combine Gabber's Absolute Purity Theorem (see \cite{Fujiwara02}) with Lem.\ \ref{lem: canonical lift of the pullback map} applied to the open subscheme given by the punctured model $\mathfrak{X}^\bullet := \mathfrak{X} \setminus Y_{\rm Sing} \hookrightarrow \mathfrak{X}$.
Note that we can not work with the pullback map to $\mathfrak{X}$ itself, since $\Delta_{\mathfrak{X}}$ is a Cartier-divisor on $\mathfrak{X} \times_\mathfrak{o} \mathfrak{X}$ if and only if $\mathfrak{X}/\mathfrak{o}$ is smooth.

\begin{prop}\label{prop: canonical lift of the cycle class to the model}
Let $X/k$ be a geometrically connected smooth projective curve over a $p$-adic field $k$ of genus $\geq 2$, $\mathfrak{X}/\mathfrak{o}$ a regular, proper, flat model s.t.\ the reduced special fibre $Y=(\mathfrak{X}\otimes_{\mathfrak{o}}\mathbb{F})_{\rm red}/\mathbb{F}$ satisfies the assumptions of Thm.\ \ref{thm: homotopy type of a curve} (iii).
Then for any $\ell \neq p$ and any section $s$ of $\pi_1(X/k)$, the induced $\ell$-adic cycle class ${\rm cl}_s$ admits a canonical lift ${\rm cl}_s^{\mathfrak{X}}$ to ${\rm H}^2(\mathfrak{X},\mathbb{Z}_\ell(1))$.
\end{prop}

\proof
First, note that $\mathfrak{X}^\bullet$ still has generic fibre $X$.
The diagonal $\Delta_{\mathfrak{X}}$ of $\mathfrak{X} \times_\mathfrak{o} \mathfrak{X}$ is Cartier outside of the collection of points $(y,y)$ for $y$ a singular point of $Y$.
It follows, that $\Delta_{\mathfrak{X}}\vert_{\mathfrak{X}^\bullet \times_\mathfrak{o} \mathfrak{X}}$ is a Cartier-divisor on $\mathfrak{X}^\bullet \times_\mathfrak{o} \mathfrak{X}$.
By construction, its Chern-class $\hat{c}_1[\mathfrak{X}^\bullet \times_\mathfrak{o} \mathfrak{X}]$ is a lift of the Chern class $\hat{c}_1[\Delta_X]$ of the diagonal.
It follows by Lem.\ \ref{lem: canonical lift of the pullback map}, that $({\rm id}_{\mathfrak{U}}^* \otimes^{\mathbb{L}} \bar{s}_\ell^*)(\hat{c}_1[\mathfrak{X}^\bullet \times_\mathfrak{o} \mathfrak{X}])$ is a canonical lift of ${\rm cl}_s = ({\rm id}_U^* \otimes^{\mathbb{L}} s^*)(\hat{c}_1[\Delta_X])$ to ${\rm H}^2(\mathfrak{X}^\bullet,\mathbb{Z}_\ell(1))$.
But $Y_{\rm Sing} \hookrightarrow \mathfrak{X}$ is a closed immersion of regular schemes of pure codimension $2$, so ${\rm H}^2(\mathfrak{X},\mathbb{Z}_\ell(1))\rightarrow {\rm H}^2(\mathfrak{X}^\bullet,\mathbb{Z}_\ell(1))$ is an isomorphism by Gabber's Absolute Purity Theorem and we get a canonical lift ${\rm cl}_s^{\mathfrak{X}}$ of ${\rm cl}_s$ to ${\rm H}^2(\mathfrak{X},\mathbb{Z}_\ell(1))$. 
\Endproof

\bibliographystyle{alpha}

%\bibliography{UniversalBibFile}%.bib}

\end{document}